\documentclass{aims} 
\usepackage{amsmath}
\usepackage{paralist}
\usepackage[misc]{ifsym}
\usepackage{epsfig} 
\usepackage{epstopdf} 
\usepackage[colorlinks=true]{hyperref}
\usepackage{adjustbox}
\usepackage{booktabs}
\usepackage{subcaption}
\hypersetup{urlcolor=blue, citecolor=red}
\allowdisplaybreaks

\def\currentvolume{X}
 \def\currentissue{X}
  \def\currentyear{200X}
   \def\currentmonth{XX}
    \def\ppages{X-XX}
     \def\DOI{10.3934/xx.xxxxxxx}

\newtheorem{theorem}{Theorem}[section]
\newtheorem{corollary}[theorem]{Corollary}

\newtheorem{lemma}[theorem]{Lemma}
\newtheorem{proposition}[theorem]{Proposition}

\theoremstyle{definition}

\newtheorem{remark}[theorem]{Remark}
\newtheorem{assumption}[theorem]{Assumption}

\title[Statistical inference for ergodic diffusion with Markovian switching]
{Statistical inference for ergodic diffusion with Markovian switching} 

\author[Yuzhong Cheng and Hiroki Masuda]{}

\subjclass{Primary~62M05, 60H10, 60J60, Secondary~60K37.}
\keywords{Asymptotic inference, diffusion process, ergodicity, Gaussian quasi-likelihood, Markovian switching.}

\thanks{$^*$Corresponding author: Hiroki Masuda}

\begin{document}
\maketitle

\centerline{\scshape
Yuzhong Cheng$^{{\href{cheng.yuzhong.129@s.kyushu-u.ac.jp}{\textrm{\Letter}}}1}$
and Hiroki Masuda$^{{\href{hmasuda@ms.u-tokyo.ac.jp}{\textrm{\Letter}}}*2}$}

\medskip

{\footnotesize
 \centerline{$^1$Joint Graduate School of Mathematics for Innovation, Kyushu University, 744 Motooka Fukuoka, Japan}
} 

\medskip

{\footnotesize
 \centerline{$^2$Graduate School of Mathematical Sciences, University of Tokyo, 3-8-1 Komaba Meguro-ku Tokyo, Japan}
}

\bigskip

 \centerline{(Communicated by Handling Editor)}


\begin{abstract}
This study explores a Gaussian quasi-likelihood approach for estimating parameters of diffusion processes with Markovian regime switching. 
Assuming the ergodicity under high-frequency sampling, we will show the asymptotic normality of the unknown parameters contained in the drift and diffusion coefficients and present a consistent explicit estimator for the generator of the Markov chain. Simulation experiments are conducted to illustrate the theoretical results obtained.
\end{abstract}


\section{Introduction}

In this paper, we consider the parameter estimation problem for Markovian switching diffusion with discrete observations.
Stochastic differential equations with Markovian switching, a subset of Markovian switching models, have attracted considerable attention across a broad range of practical applications. The Markovian switching model is a statistical model assuming that the data can be divided into different regimes or states, and transitions between these states are governed by a continuous-time finite-state Markov chain. This concept garnered attention following its introduction by Hamilton in \cite{Hamiltion1989} within the context of autoregressive models, and has since been extended to settings involving stochastic differential equations.

Switching diffusion models constitute a class of hybrid systems that integrate continuous dynamics with discrete events.
From an applied perspective, switching diffusion models are playing an increasingly pivotal role across diverse fields, including economics, financial engineering, ecology, and biological systems. For example, \cite{ZHANG2022273} applied switching diffusion models to price barrier options in financial markets, leveraging their capacity to account for structural changes in economic conditions and fluctuations in business and investment environments.
Moreover, the application of switching diffusion models has proven useful in modeling a wide range of ecological phenomena. For example, these models have been used to understand animal movement dynamics, as demonstrated in studies such as \cite{BLACKWELL199787}, \cite{Blackwell2003}, \cite{harris2013flexible}, \cite{patterson2017statistical}.

As mentioned in \cite{patterson2017statistical}, movement ecology aims to understand the reasons behind organisms' movements through space and the constraints they face during these movements. Real animal movements and behaviors are highly complex and dynamic. There are limitations to what can be inferred solely from position and sensor data. Therefore, to extract significant patterns from the data, it is often useful to assume that movement processes are driven by switches between different behavioral modes. Various modeling approaches have been developed to accommodate these different phases or modes of movement. 
Several works have illustrated the usage of (Markov) switching diffusion models in animal movement ecology. 
Initially proposed by \cite{BLACKWELL199787}, this approach utilizes the Ornstein-Uhlenbeck (OU) process with Markovian switching to characterize the movement of individual animals, further elaborated upon by subsequent researches \cite{Blackwell2003}, \cite{harris2013flexible}, and \cite{patterson2017statistical}, to mention just a few.
In this framework, animal movement is characterized by two components: a position process $U_t$ representing the location of the animal, typically described by an OU process, and a behavioral process governing movement patterns, modeled as a continuous-time Markov chain. The behavioral process  $M_t$  is defined as a continuous-time Markov chain, which models the individual's transitions between different states over time, with each state potentially representing a distinct behavior or position. The behavioral process $ M_t $  is further characterized as a continuous-time Markov chain with a finite state space $ \{1, 2, \ldots, N\} $ and a transition rate matrix $Q = (q_{ij})$, which determines the rate at which the animal switches between different behavioral states, with the entry $q_{ij}$ representing the transition rate from state $i$ to state $j$.

As stated in \cite{patterson2017statistical}, the position process $U_t$ is formally expressed as:
\begin{align*}
	dU_t = \beta_{M_t}(\gamma_{M_t} - U_t) \, dt + \sigma_{M_t} \, dw_t,
\end{align*}
where $M_t$ denotes the behavioral process and $ w_t$ is a standard Brownian motion. 
The parameters ${(\beta_1, \gamma_1, \sigma_1), \ldots, (\beta_N, \gamma_N, \sigma_N)}$ are defined, with each tuple $(\beta_i, \gamma_i, \sigma_i)$ corresponding to a specific behavioral state.
Consequently, when the animal resides in a behavioral state $i$, its movement adheres to an OU process characterized by the parameters $ \beta_{i}, \gamma_i, \sigma_i $. Furthermore, as explored in \cite{michelot2019}, this framework can be extended to incorporate velocity dynamics.
As noted in \cite{Blackwell2003} and \cite{harris2013flexible}, these studies focus on scenarios in which both behavior and position are observable, a situation that is increasingly common due to technological advancements. We also consider this.
Moreover, nonlinear diffusion models with regime switching are in principle feasible, as discussed in Section 2.4.4 of \cite{wreo25274}. As special cases in the study \cite{michelot2021varying}, nonlinear diffusion models with regime switching have been utilized to investigate specific behaviors such as elephant movement and the diving patterns of beaked whales.

Although switching diffusion models have found widespread use in various applications, the problem of parameter estimation within this framework has not received extensive attention. In a recent study \cite{Yuhang2023}, the following stochastic differential equation was considered
\begin{equation*}
	dX_{t}^{\varepsilon}=b\left(X_{t}^{\varepsilon}, \Lambda_{t}, \theta\right)dt+\varepsilon \sigma\left(X_{t}^{\varepsilon}, \Lambda_{t}\right) dB_{t}, \quad t \in[0, T].
\end{equation*}
Here,  $\varepsilon \in(0,1)$ represents the scale parameter, $B_{t}$ denotes a standard Brownian motion, and $\Lambda_t$ represents a right-continuous Markov chain taking values in a finite state space. In their work, they focused on estimating the drift parameter $\theta$ 
based on discrete observations $\{(X_{t_k}^{\varepsilon},\Lambda_{t_k})\}$ within a fixed time interval $[0, T]$, where $\{t_k=k\Delta, k=1,2,...,n\}$ as $\varepsilon \to 0$ and $\Delta \to 0$. They derived the consistency and asymptotic distribution of the least-squares estimator.

In this paper, we consider an ergodic counterpart of the aforementioned study. Let $Z_t$ be a finite state continuous-time Markov chain with generator $Q$ and $w_t$ represents a standard Wiener process independent with $Z_t$. We consider the parameter estimation problem for both drift and diffusion parameters within the context of Markovian switching dependence, governed by the following one-dimensional stochastic differential equation
\begin{equation}
	dX_t = b(X_t,Z_t, \alpha_{Z_t})dt + \sigma(X_t,Z_t, \gamma_{Z_t})dw_t.
	\label{eq:SDEmodel}
\end{equation}
Our analysis is based on discrete observations $\{(X_{t_0},Z_{t_0}),(X_{t_1},Z_{t_1}),...,(X_{t_n}Z_{t_n})\}$ where $j=0,1,..,n$, $t_j = jh$, with $h=h_n$ denoting the sampling stepsize, under high-frequency sampling conditions:
\begin{equation}\nonumber
	\text{$T_n := nh \to \infty$\quad and \quad $nh^2 \to 0$\quad as $n \to \infty$.}    
\end{equation}
Assuming the ergodicity of the solution process, we will employ a Gaussian quasi-likelihood inference, as presented in \cite{Kes97}, and prove the consistency and asymptotic normality of the associated estimator.

Furthermore, we also propose a consistent estimator for generator $Q$ based solely on observation $(Z_{t_j})_{j=0}^{n}$. 
There is extensive literature on the estimation of $Q$ in various settings and perspectives. 
In \cite{bladt2005statistical}, the authors examined the properties of the maximum likelihood estimator (MLE) of the generator $Q$ and established the consistency of the MLE for ergodic Markov chains based on discrete observations under a fixed small time step size. Building on this result, they also provided an EM algorithm and an MCMC approach to estimate the generator $Q$. \cite{dos2018robust} improved the EM algorithm discussed in \cite{bladt2005statistical} by offering directly computable closed-form expressions for the quantities appearing in the EM algorithm.
As alternative approaches to this estimation problem, \cite{metzner2007generator} provides a summary and performance comparison of several existing methods for estimating the generator of a continuous-time Markov chain.
However, unlike other estimation methods that rely on numerical algorithms, the quasi-likelihood function approach to estimate the generator $Q$ does not seem to have been investigated as yet, while it is very natural in the high-frequency-sampling scenario. 
Our quasi-likelihood approach to estimate $Q$ provides us with an easily computable estimator; see \eqref{eq:quasiQ} for the definition.

Finally, we note that the theoretical results in this paper can be extended to the multidimensional $X$ without any essential change.

The paper is organized as follows. Section \ref{sec:setup} outlines the model setting for diffusions with Markovian switching, the high-frequency observation framework, and key assumptions. Section \ref{sec:gqli} presents the main results, namely the Gaussian quasi-likelihood inference for our model. 
Section \ref{sec:gen.Q} constructs an explicit, easily computable estimator for the generator of the Markov chain and proves its consistency under high-frequency observation. Section \ref{sec:sim} presents numerical experiments to evaluate the performance of our estimators. Finally, all proofs are provided in Section \ref{sec:proofs}.

\section{Model setup}
\label{sec:setup}

Let $(\Omega,\mathcal{F},\{\mathcal{F}_t\}_{t \geq 0},\mathbb{P})$ be a complete filtered probability space satisfying the usual conditions, on which all the random quantities appearing below are defined. We consider the one-dimensional switching diffusion model given by \eqref{eq:SDEmodel}, where the ingredients are described below.
\begin{itemize}
	\item $w=(w_t)$ is a one-dimensional real-valued Wiener process.
	\item $Z=(Z_t)$ is a homogeneous continuous-time Markov chain with finite state space $S := \{1,2,..., N\}$, where $N$ is the number of the states in this Markov chain, which is a finite positive integer, and the generator $Q = (q_{ij})$ of this Markov chain $Z$ is a $N$ by $N$ matrix  characterized by the Kolmogorov backward equation
	\begin{align*}
		\frac{dP(t)}{dt} = QP(t), \,\qquad P(0) = I_N,
	\end{align*} where $I_N$ is the $N\times N$-identity matrix. The elements of matrix $Q$ are given by
	\begin{equation}
		\begin{array}{c}
			q_{ij}
			=\left\{\begin{array}{ll}
				\lim_{h \to 0} \frac{1}{h}\mathbb{P}(Z_{t+h}=j \mid Z_t=i) , & \text { if } i \neq j, \\[2mm]
				\lim_{h \to 0} \frac{1}{h}\{\mathbb{P}(Z_{t+h}=j \mid Z_t=i)-1\}, & \text { if } i = j,
			\end{array}\right.
		\end{array}
		\label{CTMCgenerator}
	\end{equation}
	with the property 
	\begin{equation}\nonumber
		q_{ii} = - \sum_{j \in S, j \neq i} q_{ij}.
	\end{equation}
	We note that the Markov chain $Z$ can be represented as a stochastic integral with respect to a Poisson random measure (see \cite[Section 1.7]{mao2006stochastic} and \cite{shao2018Euler}). 
	Suppose for a moment that $q_{ij}>0$ for $i\ne j$; this will be assumed later.
	For $i, j \in S$  with  $j \neq i$, let  $\Delta_{i j}$  be consecutive, left closed and right open intervals of the real line, each having length  $q_{i j}$, defined as follows:
	\begin{align*}
		&\Delta_{12} = [0,q_{12}), \: \Delta_{13} = [q_{12},q_{12}+q_{13}) \: 
		\\
		&... \:  
		\\
		& \Delta_{1N}=\left[\sum_{j=2}^{N-1}q_{1j},\sum_{j=2}^{N}q_{1j}\right),
		\\
		&\Delta_{21}=\left[\sum_{j=2}^{N}q_{1j},\sum_{j=2}^{N}q_{1j}+q_{21}\right), \: \Delta_{23}=\left[\sum_{j=2}^{N}q_{1j}+q_{21},\sum_{j=2}^{N}q_{1j}+q_{21} + q_{23}\right),
		\\
		&...
		\\
		&\Delta_{2N}=\left[\sum_{j=2}^{N}q_{1j}+\sum_{j=1,j\neq2}^{N-1}q_{2j},\sum_{j=2}^{N}q_{1j}+\sum_{j=1,j\neq2}^{N}q_{2j}\right),
	\end{align*} and so on. Let $I_A$ be the indicator function of a set $A$.
	Define a function  $h: S \times \mathbb{R} \mapsto \mathbb{R} $ by
	\begin{equation*}
		h(i , y)=\sum_{j \ne i}(j-i) I_{\left\{y \in \Delta_{i j}\right\}} .
	\end{equation*}
	This implies that for each   $i \in S$ , if $y \in \Delta_{i j}$ , then $h( i, y)=j-i $; otherwise,  $h( i, y)=0 $. Then, as in \cite{shao2018Euler}
	\begin{equation}
		dZ_t=\int_{\mathbb{R}} h\left(Z_{t-}, y\right) \Pi(dt, dy),
		\nonumber
	\end{equation}
	where  $\Pi(dt, dy) $ is a Poisson random measure with intensity  $d t \times m(dy)$ , and  $m(\cdot) $ is the Lebesgue measure on $ \mathbb{R}$ . The Poisson random measure  $\Pi(\cdot, \cdot)$ here is independent of the Wiener process $w$ in equation \eqref{eq:SDEmodel}.
	
	\item For each $i \in S$, $\alpha_i \in \Theta_{\alpha_i} \subset \mathbb{R}, \gamma_i \in \Theta_{\gamma_i} \subset \mathbb{R}$ are parameters depend on the state of the Markov chain $Z$, where $\Theta_{\alpha_i}$, $\Theta_{\gamma_i}$ are compact convex sets. Define $\Theta_{\alpha}:=\prod_i\Theta_{\alpha_i}$, $\Theta_{\gamma}:=\prod_i\Theta_{\gamma_i}$ and the whole parameter space $\Theta:= \Theta_{\alpha} \times \Theta_{\gamma}$.
	\item $b(x,i,\alpha_i)$ and $\sigma(x,i,\gamma_i)$ are real-valued measurable functions.
	\item $w$, $Z$, and the initial value $X_0$ are mutually independent.
\end{itemize}

Let $\alpha = (\alpha_1,...,\alpha_N)$, $\gamma = (\gamma_1,...,\gamma_N)$ and $\theta := (\alpha,\gamma)$. We estimate the parameter $\theta=(\alpha,\gamma)$ with discrete observations under high-frequency sampling. Let $j=0,1,..,n$, $t_j = jh$, with $h=h_n$ denoting the sampling stepsize. We observe a sample $\{(X_{t_0},Z_{t_0}),(X_{t_1},Z_{t_1}),...,(X_{t_n},Z_{t_n})\}$ with 
\begin{equation}\nonumber
	\text{$T_n = nh \to \infty$ \quad and \quad $nh^2 \to 0$\quad  as $n \to \infty$.}
\end{equation}
We denote the true value of parameter $\theta$ be $\theta^{\star}$, which is assumed to lie in the interior of $\Theta$. We use $\mathbb{P}_{\theta}$ to represent the distribution of the process $(X, \alpha)$ with the assigned parameter value $\theta$, and $\mathbb{E}_{\theta}$ denotes the expectation with respect to the distribution $\mathbb{P}_{\theta}$.  The symbol $\xrightarrow{p}$ denotes convergence in probability with respect to $\mathbb{P}_{\theta^{\star}}$.

Define the differential operator $\partial_x:= \frac{\partial}{\partial x}$ as the derivative with respect to variable $x$, and let $\partial_x^k$ denote the derivative with respect to $x$ up to order $k$. 
We introduce the following assumption to ensure the existence and uniqueness of the solution to \eqref{eq:SDEmodel}.

\begin{assumption} We make the following assumptions on drift and diffusion coefficients:
	\label{ass:smoothness}
	\begin{enumerate}
		\item There exists a constant $C >0 $ such that 
		\begin{align}
			&|b(x,i,\alpha_i)-b(y,i,\alpha_i)| + |\sigma(x,i,\gamma_i)-\sigma(y,i,\gamma_i)| \leq C|x-y|,
			\\
			&|b(x,i,\alpha_i)|^2 +|\sigma(x,i,\gamma_i)|^2 \leq C(1+|x|^2)
		\end{align} for any $x,y \in \mathbb{R}$ and each $i \in S$.
		\item  The coefficients $b(x,i,\alpha_i), \sigma(x,i,\gamma_i))$ are twice continuously differentiable with respect to the first variable and three times continuously differentiable with respect to the third variable, with $\sigma(x,i,\gamma_i)>0$. There exists a nonnegative constant $C$ satisfying that 
		\begin{align}
			&\max_{i\le N}\sup _{(x,\alpha_i,\gamma_i) \in \mathbb{R} \times \Theta_{\alpha_i} \times \Theta_{\gamma_i}} \frac{1}{1+|x|^C} \left(\left|\partial_{\alpha_i}^{k} \partial_{x}^{l} b(x, i, \alpha_i)\right| + \left|\partial_{\gamma_i}^{k} \partial_{x}^{l} \sigma(x, i, \gamma_i)\right|
			+ \sigma^{-1}(x, i, \gamma_i)\right) 
			\nonumber\\
			&< \infty,
		\end{align}  where $k \in \{0,1,2,3\}$ and $l \in \{0,1,2\}$.
	\end{enumerate}  
\end{assumption}
Under Assumption \ref{ass:smoothness}, there exists a unique solution to \eqref{eq:SDEmodel}, see \cite[Section 3.2]{mao2006stochastic} and \cite[Chapter 2]{yin2009hybrid}.

The generator $Q$ is called \textit{irreducible} if the system of $N$ equations $\nu Q = 0$ subject to $\sum_{i=1}^{N}\nu_i=1$ has a unique solution $\nu=(\nu_1,\dots,\nu_N)$; see \cite[Definition A.7]{yin2009hybrid}.
We also give some assumptions on the ergodicity and boundedness of the moments of the solution process.
\begin{assumption}
	We make the following assumptions:
	\begin{enumerate}
		\item  The generator $Q$ is irreducible.
		\item 
		The process $(X,Z)$ admits a unique invariant probability measure $\nu_{\theta^{\star}}(dx,i)$ and the distribution of $(X_0,Z_0)$ is $\nu_{\theta^{\star}}$.
		\item For all $p \geq 0$, $\mathbb{E}_{\theta^{\star}}(|X_0|^p) < \infty$.
	\end{enumerate}
	\label{Ass:finite moment}
\end{assumption}
Recall that $Z_t$ has finite states. Under Assumption \ref{Ass:finite moment}, there is a unique stationary distribution $\pi = (\pi_1,...,\pi_N)$ of $Z_t$ where $\pi Q = 0$ subject to $\sum_{i=1}^{N}\pi_i=1$, and $Z_0$ have distribution $\pi$.
The criteria for the existence and uniqueness of the invariant probability measure can be found in, for example, \cite{bao2016approximation}, \cite{YUAN2003277}, and \cite[Chapter 4]{yin2009hybrid}.
Under Assumption \ref{Ass:finite moment}, \cite[Theorem 4.4]{yin2009hybrid} gives the ergodic theorem: 
for $\rho \in \{\alpha, \gamma\}$,
\begin{equation}
	\frac{1}{T}\int_{0}^{T} f(X_s,Z_s,\rho_{Z_s}) ds \xrightarrow{p} \sum_{i =1}^{N}\int_{\mathbb{R}}f(x,i,\rho_i)\nu_{\theta^{\star}}(dx,i)
	\label{eq:c-ergodic}
\end{equation}
for any real-valued Borel measurable function $f(x,i,\rho_i)$ such that
\begin{equation*}
	\sum_{i =1}^{N}\int_{\mathbb{R}}|f(x,i,\rho_i)|\nu_{\theta^{\star}}(dx,i) < \infty.
\end{equation*}

We also make the following assumptions on the identifiability of the parameter.

\begin{assumption}
	We assume that if $b(x,i,\alpha_i)=b(x,i,\alpha_i^{\star})$ for almost sure $x$ and each $i\in S$, then $\alpha=\alpha^{\star}$. Moreover, if $\sigma(x,i,\gamma_i)^2=\sigma(x,i,\gamma_i^{\star})^2$ for almost sure $x$ and each $i\in S$, then $\gamma=\gamma^{\star}$.
	\label{ass:identi}
\end{assumption}

Let $J_1$ denote the first jump time of Markov chain $Z_t$, Following \cite{liu2015central}, we define 
\begin{equation*}
	\tau_i := \inf \{t>J_1|Z_t = i\}
\end{equation*} to be the \textit{first return time} of $Z_t$ to state $i \in S$.
Let $\mathbb{P}_Q^h(i, k)$ denote the transition probability of $Z$ from state $i$ to state $k$ during time $h$.
Since our Markov chain $Z_t$ is irreducible with finite states, it is automatically positive recurrent. We say $Z_t$ is \textit{exponentially ergodic} if
\begin{equation*}
	e^{ct} \sum_{k \in S} \left|\mathbb{P}_Q^t(i, k )- \pi_k\right|  \to 0
\end{equation*} as $t \to \infty$ for some $c >0 $ and for any $i \in S$.

\section{Gaussian quasi-likelihood inference}
\label{sec:gqli}

Let $\phi(\cdot;\mu,\sigma^2)$ denote the Gaussian density function with mean $\mu$ and variance $\sigma^2$, and let
\begin{align}
	& \textbf{e}_{j-1}:=(I_{\{Z_{t_{j-1}}=1\}},...,I_{\{Z_{t_{j-1}}=N\}})
	\notag\\
	&\mu_{j-1}(\alpha) := X_{t_{j-1}} + b(X_{t_{j-1}},Z_{t_{j-1}},\alpha \cdot \textbf{e}_{j-1})h =X_{t_{j-1}} + b(X_{t_{j-1}},Z_{t_{j-1}},\alpha_{Z_{t_{j-1}}})h,
	\notag\\
	&\sigma_{j-1}(\gamma) := \sigma(X_{t_{j-1}},Z_{t_{j-1}},\gamma \cdot \textbf{e}_{j-1}) =\sigma(X_{t_{j-1}},Z_{t_{j-1}},\gamma_{Z_{t_{j-1}}}).
	\notag
\end{align}
Our construction of the quasi-likelihood goes as follows:
since $(X,Z)$ is Markov, we may formally write the log-likelihood as $\theta \mapsto \sum_{j=1}^{n} \log f_{n,\theta}(X_{t_{j}},Z_{t_{j}}|X_{t_{j-1}},Z_{t_{j-1}})$ for some conditional density $f_{n,\theta}(x,z|x',z')$. Since $Z_{t_{j-1}}=Z_{t_j}$ with high probability in small time, it is natural to consider the following approximate version, say $\theta \mapsto \sum_{j=1}^{n} \log f_{n,\theta}(X_{t_{j}}|X_{t_{j-1}},Z_{t_{j-1}})$. This observation leads to the Gaussian quasi-likelihood function
\begin{align}
	\mathbb{H}_n(\theta) &:= \sum_{j=1}^{n} \log \phi(X_{t_j};\mu_{j-1}(\alpha),h\sigma_{j-1}^2(\gamma))
	\notag\\
	&=C_n - \frac{1}{2}\sum_{j=1}^{n} \left( \log \sigma_{j-1}^2(\gamma) + \frac{\left(X_{t_j}-\mu_{j-1}(\alpha)\right)^2}{h\sigma_{j-1}^2(\gamma)} \right)
	\nonumber\\
	&= C_n - \frac12 \sum_{j=1}^{n} \sum_{i\in S} \left( \log \sigma_{j-1}^2(\gamma) + \frac{\left(X_{t_j}-\mu_{j-1}(\alpha)\right)^2}{h\sigma_{j-1}^2(\gamma)} \right)I_{\{Z_{t_{j-1}}=i\}},
	\label{eq:GQLF}
\end{align}
where $C_n$ is a constant that only depends on $n$.
Then Gaussian quasi maximum likelihood estimator (GQMLE) $\hat{\theta}_n$ is defined by
\begin{equation}
	\hat{\theta}_n \in \mathop{\rm argmax}_{\theta \in \Theta} \mathbb{H}_n(\theta).
	\label{eq:GQMLE}
\end{equation}
We use the following notation for estimators above $\hat{\theta}_n = (\hat{\alpha}_n,\hat{\gamma}_n)$, and $\hat{\alpha}_n =(\hat{\alpha}_{1,n},...,\hat{\alpha}_{N,n})$, $\hat{\gamma}_n = (\hat{\gamma}_{1,n},...,\hat{\gamma}_{N,n})$.

\begin{remark}
	If we know that some components of $\alpha$ and/or $\gamma$ are the same, then we could incorporate it by merging the summation over $i\in S$ in \eqref{eq:GQLF}; 
	for example, we could allow $\alpha_1^\star=\alpha_2^\star$ for $Z_t\in\{1,2\}$.
	Also, we could handle the model
	\begin{equation}\nonumber
		dX_t = b(X_t,Z_t, \alpha)dt + \sigma(X_t,Z_t, \gamma)dw_t,
	\end{equation}
	where the parameters $\alpha$ and $\gamma$ are constant so that they do not varying according as $Z$.
\end{remark}

Let 
\begin{gather*}
	D_n := \operatorname{diag}(\underbrace{\sqrt{nh},\ldots, \sqrt{nh}}_{N \text{ times}}, \underbrace{\sqrt{n},\ldots, \sqrt{n}}_{N \text{ times}}),
	\\
	\mathcal{I}(\theta^{\star}):= \operatorname{diag}\left(G_{1,1}(\theta^{\star}),\ldots, G_{1,N}(\theta^{\star}), G_{2,1}(\theta^{\star}),\ldots, G_{2,N}(\theta^{\star})\right),
\end{gather*}
where for $i \in S$
\begin{align*}
	&G_{1,i}(\theta^{\star}) := \int_{\mathbb{R}}\left(\frac{\left(\partial_{\alpha_i}b(x,i,\alpha_i^{\star})\right)^2}{\sigma^2(x,i,\gamma_i^{\star})}\right) \nu_{\theta^{\star}}(dx,i),
	\\
	&G_{2,i}(\theta^{\star}) := \int_{\mathbb{R}}\left(\frac{\left(\partial_{\gamma_i}\sigma^2(x,i,\gamma_i^{\star})\right)^2}{2\sigma^4(x,i,\gamma_i^{\star})}\right) \nu_{\theta^{\star}}(dx,i).
\end{align*}
We now state the main result of this paper.
\begin{theorem}
	\label{thm:AN}
	Under Assumptions \ref{ass:smoothness} to \ref{ass:identi}, we have
	\begin{align}
		D_n(\hat{\theta}_n-\theta^{\star}) \xrightarrow{\mathcal{L}} N(0,\mathcal{I}(\theta^{\star})^{-1}).
	\end{align}
\end{theorem}

Based on the expression of $\mathcal{I}(\theta^{\star})$, we can construct an estimator $\hat{\mathcal{I}}_n$ for $\mathcal{I}(\theta^{\star})$:
\begin{align*}
	\hat{\mathcal{I}}_n := \operatorname{diag}\left(\hat{G}_{1,1}^{(n)},\ldots, \hat{G}_{1,N}^{(n)}, \hat{G}_{2,1}^{(n)},\ldots, \hat{G}_{2,N}^{(n)}\right),
\end{align*}
where
\begin{align*}
	&\hat{G}_{1,i}^{(n)} := \frac{1}{n} \sum_{j=1}^{n} \frac{\left(\partial_{\alpha_i}b(X_{t_{j-1}},Z_{t_{j-1}},\hat{\alpha}_{i,n})\right)^2}{\sigma^2(X_{t_{j-1}},Z_{t_{j-1}},\hat{\gamma}_{i,n})} I_{\{Z_{t_{j-1}}=i\}},
	\\
	&\hat{G}_{2,i}^{(n)} := \frac{1}{n} \sum_{j=1}^{n}\frac{\left(\partial_{\gamma_i}\sigma^2(X_{t_{j-1}},Z_{t_{j-1}},\hat{\gamma}_{i,n})\right)^2}{2\sigma^4(X_{t_{j-1}},Z_{t_{j-1}},\hat{\gamma}_{i,n})} I_{\{Z_{t_{j-1}}=i\}},
\end{align*}
for $i \in S$. The following result is a consequence of Theorem \ref{thm:AN} and Lemma \ref{lem:disergodic}.
\begin{corollary}
	Under Assumptions \ref{ass:smoothness} to \ref{ass:identi}, the consistency
	$\hat{\mathcal{I}}_n \xrightarrow{p} \mathcal{I}(\theta^{\star})$ holds.
	\label{cor:consistencyI}
\end{corollary}

Thus, we have
\begin{equation}
	\hat{\mathcal{I}}_n^{1/2} D_n(\hat{\theta}_n-\theta^{\star}) \xrightarrow{\mathcal{L}} N(0,I_{2N}),
\end{equation}
based on which we can consider testing the equality of some components of the parameters; for example, we can consider testing the hypothesis $\alpha_1=\alpha_2$ and so on, more generally the linear hypothesis $A\theta=0$ for a given matrix $A$.

\subsection{An example: Ornstein-Uhlenbeck switching diffusion}
\label{ex:ou}

We consider the following one-dimensional real-valued Ornstein-Uhlenbeck switching diffusion
\begin{align}
	dX_t = - \beta_{Z_t}X_tdt + \sigma_{Z_t}dw_t,
	\label{eq:OUp}
\end{align}
where $Z_t$ is continuous-time Markov chain and state space $S=\{1,2,...,N\}$, For each $i \in S$, $\beta_i \in \Theta_{\beta_i} \subset (0,\infty), \sigma_i \in \Theta_{\sigma_i} \subset (0,\infty)$ are parameters depend on the state of the Markov chain $Z_t$, where $\Theta_{\beta_i}$, $\Theta_{\sigma_i}$ are compact convex set.
Let $\beta = (\beta_1,...,\beta_N)$, $\sigma = (\sigma_1,...,\sigma_N)$ and $\theta := (\beta,\sigma)$ and true parameter value be $\theta^{\star}$. 
The regularity assumptions in Assumption \ref{ass:smoothness} are easily satisfied.
Direct computations give, for $i \in S$,
\begin{align}
	&\hat{\beta}_{i,n} = \frac{\sum_{j=1}^{n} \left(X_{t_{j-1}}^2-X_{t_j}X_{t_{j-1}}\right)I_{\{Z_{t_{j-1}}=i\}}}{h\sum_{j=1}^{n} X_{t_{j-1}}^2 I_{\{Z_{t_{j-1}}=i\}}},
	\label{estor:ou_1}\\
	&\hat{\sigma}^2_{i,n} = \frac{\sum_{j=1}^{n}\left(X_{t_j}-X_{t_{j-1}}+\hat{\beta}_{i,n}X_{t_{j-1}}h\right)^2I_{\{Z_{t_{j-1}}=i\}}}{h \sum_{j=1}^{n} I_{\{Z_{t_{j-1}}=i\}}}.
	\label{estor:ou_2}
\end{align}
According to Theorem \ref{thm:AN}, the asymptotic covariance matrix $\mathcal{I}(\theta^{\star})^{-1}$ is given through integrals with respect to the invariant probability measure $\nu_{\theta^{\star}}$.
The explicit form of invariant probability measure for switching diffusion is usually difficult to obtain due to the presence of switching, as mentioned in \cite[Chapter 1]{yin2009hybrid}. There are several studies on numerical approximation of invariant probability measure, for example, \cite{bao2016approximation}, \cite{yin2005numerical} and \cite{yin2009hybrid}.

For the Ornstein-Uhlenbeck switching diffusion \eqref{eq:OUp}, the explicit form of invariant probability measure when state number $N=2$ is known. Let $N=2$ in the sequel and let $\pi=(\pi_1,\pi_2)$ be the stationary distribution of $Z_t$, which is ensured when generator $Q$ of $Z_t$ is irreducible. Then \cite[Theorem 1.1]{Zhang2017} gives the conditions for the existence and uniqueness of invariant probability measure for equation \eqref{eq:OUp}. 
Under the condition
\begin{equation}
	\frac{\sigma_1}{2\beta_1} = \frac{\sigma_2}{2\beta_2},
	\label{con:2ou}
\end{equation}
Theorem 3.1 in the same paper gives the explicit form of the invariant probability measure through a Fourier transform:
\begin{align*}
	\nu_{\theta}(dx,k) &= \frac{1}{2\pi} \int_{-\infty}^{+\infty}\pi_k \exp{\left(-\frac{\sigma_k^2}{2\beta_k}\xi^2-i\xi x\right)} d\xi dx
	\\
	&= \frac{\pi_k}{\sqrt{2\pi(\sigma_k^2/\beta_k)}} \exp\left(-\frac{\beta_kx^2}{2\sigma_k^2}\right)
\end{align*}
for $k \in \{1,2\}$. Then, the asymptotic covariance $\mathcal{I}(\theta^{\star})^{-1}$ is given by
\begin{align}
	\mathcal{I}(\theta^{\star})^{-1} = \operatorname{diag}\left(G_{1,1}(\theta^{\star}), G_{1,2}(\theta^{\star}), G_{2,1}(\theta^{\star}), G_{2,2}(\theta^{\star})\right)^{-1},
	\label{eq:OUav}
\end{align}
where the entries are explicitly given by
\begin{align*}
	&G_{1,k}(\theta^{\star}) = \frac{\pi_k}{\sqrt{2\pi((\sigma_k^{\star})^2/\beta_k^{\star})}}\int_{\mathbb{R}}\frac{x^2}{(\sigma_{k}^{\star})^2} \exp\left(-\frac{\beta_k^{\star} x^2}{2(\sigma_k^{\star})^2}\right) dx = \frac{\pi_k}{\beta_k^{\star}},
	\\
	&G_{2,k}(\theta^{\star}) = \frac{\pi_k}{\sqrt{2\pi((\sigma_k^{\star})^2/\beta_k^{\star})}}\int_{\mathbb{R}}\frac{2}{(\sigma_k^{\star})^2} \exp\left(-\frac{\beta_k^{\star}x^2}{2(\sigma_k^{\star})^2}\right) dx = \frac{2\pi_k}{(\sigma_k^{\star})^2}
\end{align*} for $k \in \{1,2\}$.

For further analysis, we provide an estimator for $\pi$ and establish its consistency and asymptotic normality in this example.
Under Assumption \ref{Ass:finite moment}, due to Lemma \ref{lem:ZIDL} and the ergodic theorem of continuous-time Markov chains (see, for example,  \cite[Theorem 3.8.1]{norris1998markov}), the estimators
\begin{equation}
	\hat{\pi}_k := \frac1n \sum_{j=1}^{n} I_{\{Z_{t_j}=k\}},\qquad k=1,2,
	\label{eq:esti_pi}
\end{equation}
are consistent: $(\hat{\pi}_1,\hat{\pi}_2) \xrightarrow{p} (\pi_1,\pi_2)$.
Hence, we can construct an approximate confidence intervals of $\beta^\star_k$ and $\sigma^\star_k$ through the Studentizations:
\begin{align}
	\sqrt{\frac{\hat{\pi}_k}{\hat{\beta}_k}}\sqrt{T_n}(\hat{\beta}_k - \beta_k^\star) \xrightarrow{\mathcal{L}} N(0,1), \label{hm:ex.eq-1} \\
	\sqrt{\frac{2\hat{\pi}_k}{(\hat{\sigma}_k)^2}}\sqrt{n}(\hat{\sigma}_k - \sigma_k^\star) \xrightarrow{\mathcal{L}} N(0,1).
	\label{hm:ex.eq-2}
\end{align}
Additionally, recall that  $\tau_{i}$ is the first return time to state $i$. With a stronger assumption on the ergodicity of $Z_t$, we can state the following proposition:
\begin{proposition} \label{prop:ex1}
	Assume the continuous-time Markov chain $Z_t$ is exponentially ergodic. For any $i_0 \in \{1,2\}$, $k \in \{1,2\}$ the following holds:
	\begin{align}
		\sqrt{nh}(\hat{\pi}_k-\pi_k) \xrightarrow{\mathcal{L}} N(0,V_{\pi_k}),
	\end{align}
	where
	\begin{align*}
		&V_{\pi_k}= 2 \sum_{i \in \{1,2\}} (I_{\{i=k\}}-\pi_k)F_i\pi_i ,
		\\
		&F_i = \mathbb{E}_{Q^{\star}}\left(\int_{0}^{\tau_{i_0}}(I_{\{Z_s=k\}}-\pi_k)ds \,\middle|\, Z_0 = i\right).
	\end{align*}
\end{proposition}

\begin{remark}
	When $Z_t$ is a general continuous-time Markov chain as in \eqref{eq:SDEmodel}, Proposition \ref{prop:ex1} remains valid for the estimators
	\begin{equation*}
		\hat{\pi}_k := \frac1n \sum_{j=1}^{n} I_{\{Z_{t_j}=k\}},\qquad k \in S.
	\end{equation*}
\end{remark}

\section{Estimation of generator $Q$}
\label{sec:gen.Q}

In this section, we propose a quasi-likelihood function approach to estimate $Q$ from observations of $\{Z_{t_0},Z_{t_1},...,Z_{t_n}\}$ in our high-frequency setting: $h \to 0$, $T_n = nh \to \infty$ and $nh^2 \to 0$ as $n \to \infty$.
As we know, the diagonal elements in $Q$ are fully determined by the off-diagonal elements as follows:
\begin{equation*}
	q_{ii} = - \sum_{j \in S, j \neq i} q_{ij},  \text{ for } i \in S.
\end{equation*}
Let $Q'$ denote the off-diagonal elements of the matrix $Q$ and $Q^{\star}=(q^\star_{ik})_{i,k}$ denote the true value of $Q$. The parameter space is defined as
\begin{equation}\nonumber
	\Theta_{Q} = \left\{ Q' = (q_{ij})_{i \neq j} \in \mathbb{R}^{N\times (N-1)}
	\middle|\,q_{ij} >0 \text{ for all }i,j \in S \text{ and } i \neq j \right\}.
\end{equation} 
We implicitly suppose that $\Theta_{Q}$ is an open set and that $(q^\star_{ik})_{i\ne k}\subset\Theta_Q$. 
Note that the dimension of $\Theta_{Q}$ is $N(N-1)$.

The likelihood function is given by
\begin{align*}
	L_n(Q') = \prod_{j=1}^{n} \mathbb{P}_{Q'}(Z_{t_j} \mid Z_{t_{j-1}})
	=\prod_{i=1}^{N}\prod_{k=1}^{N} \mathbb{P}_{Q'}^h(i, k)^{K_{ik}(n)},
\end{align*}
where $\mathbb{P}_{Q'}$ represents the probability measure associated with $Q$, given that $Q$ is entirely characterized by $Q'$.
Moreover, $K_{ik}(n)$ denotes the total number of transitions from state $i$ to state $k$ in the embedded Markov chain $\{Z_{t_0}, \ldots, Z_{t_n}\}$:
\begin{equation}
	K_{ik}(n) = \sum_{j=1}^{n} I_{\{Z_{t_{j-1}}=i\}}I_{\{Z_{t_j}=k\}}.
	\nonumber
\end{equation}
See, for example, \cite{bladt2005statistical}.
Note that $Z_t$ is assumed to be homogeneous as discussed in Section \ref{sec:setup}, so we have $\mathbb{P}_{Q'}(Z_{t_j}=k \mid Z_{t_{j-1}}=i) = \mathbb{P}_{Q'}^h(i, k)$ for all $j$.
Since $\mathbb{P}_{Q'}^h(i, k)$ is typically unknown, the maximum likelihood estimator of $L_n(Q')$ cannot be computed explicitly. 

Here, we propose a quasi-likelihood function that provides an easily computable estimator. 
By \eqref{CTMCgenerator}, we know
\begin{align*}
	\begin{array}{c}
		\mathbb{P}_{Q'}^h(i, k)
		=\left\{\begin{array}{ll}
			q_{ik}h +o(h) , & \text { if } i \neq k, \\
			1 -\left(\sum_{l \neq i} q_{il}\right)h +o(h), & \text { if } i = k.
		\end{array}\right.
	\end{array}
\end{align*}
Ignoring the small order term $o(h)$ (for $h\to 0$), we replace $\mathbb{P}_{Q'}^h(i, k)$ by 
\begin{align}
	F(q_{ik})=F_h(q_{ik};i,k) &:= 
	\left(1 -\sum_{l \neq i} q_{il}h\right)I_{\{i=k\}} + q_{ik}h I_{\{i \neq k\}}
	\nonumber
\end{align}
in the function $L_n(Q')$, and then take the logarithm to obtain the quasi-likelihood function:
\begin{align}
	\mathbb{M}_n(Q') :=
	\sum_{i=1}^{N} \sum_{k=1}^{N} K_{ik}(n) \log \left(F(q_{ik})\right).
	\label{eq:quasiQ}
\end{align}
We define the quasi-maximum likelihood estimator $\hat{Q}_n = (\hat{q}^{(n)}_{ij})$ for $Q$ by
\begin{equation}
	(\hat{q}^{(n)}_{ij})_{i \neq j} \in \mathop{\rm argmax}_{(q_{ij})_{i \neq j} \in \bar{\Theta}_{Q}} \mathbb{M}_n(Q), \quad
	\hat{q}^{(n)}_{ii} = - \sum_{l \neq i} \hat{q}^{(n)}_{il}.
	\label{eq:estQ}
\end{equation}
Then the estimators can be computed explicitly as
\begin{align}
	\begin{array}{c}
		\hat{q}^{(n)}_{ik}
		=\left\{
		\begin{array}{ll}
			\left(h \sum_{l=1}^{N}K_{il}(n)\right)^{-1}K_{ik}(n) & (i \neq k) \\
			- \sum_{l \neq i} \hat{q}^{(n)}_{il} = \left(h \sum_{l=1}^{N}K_{il}(n)\right)^{-1}K_{ik}(n)- h^{-1}
			& (i = k)
		\end{array}
		\right.
	\end{array}
	\label{eq:estQn}
\end{align}
for $i,k \in S$.

The estimator $\hat{q}^{(n)}_{ik}$ is consistent:
\begin{theorem}
	Suppose the stability conditions on $Z$ given by Assumption \ref{Ass:finite moment}:
	\begin{enumerate}
		\item  The generator $Q$ is irreducible;
		\item The process $Z$ admits a unique invariant discrete probability measure $\int_{\mathbb{R}}\nu_{\theta^\star}(dx,\cdot)$ on $S$, and $Z$ is strictly stationary.
	\end{enumerate}    
	Then, the consistency $\hat{q}^{(n)}_{ik} \xrightarrow{p} q^{\star}_{ik}$ holds for all $i,k \in S$.
	\label{thm:consisQ}
\end{theorem}

\section{Simulation study}
\label{sec:sim}

In this section, we conduct simulation studies to test the proposed GQMLE on equation \eqref{eq:OUp}. The parameter $\theta=(\beta_1,\beta_2,\sigma_1,\sigma_2)$ are set to be $\beta_1=1$, $\beta_2 = 2$, $\sigma_1=0.1$, and $\sigma_2 =0.2$.
and the generator
\begin{align*}
	Q =
	\begin{pmatrix}
		-0.01 & 0.01 \\
		0.01 & -0.01 \\
	\end{pmatrix}.
\end{align*}
We see that the above parameters satisfy the condition in \eqref{con:2ou}, so the process admits a unique explicit invariant probability measure, as presented in Example \ref{ex:ou}.
The stationary distribution $\pi$ for $Z$ is obtained as $\pi=(\pi_1,\pi_2) = (0.5, 0.5)$ by solving the equation $\pi Q = 0$ with the condition $\pi_1 + \pi_2 = 1$.

\begin{figure}[htp]
	\centering
	\includegraphics[width=0.8\textwidth]{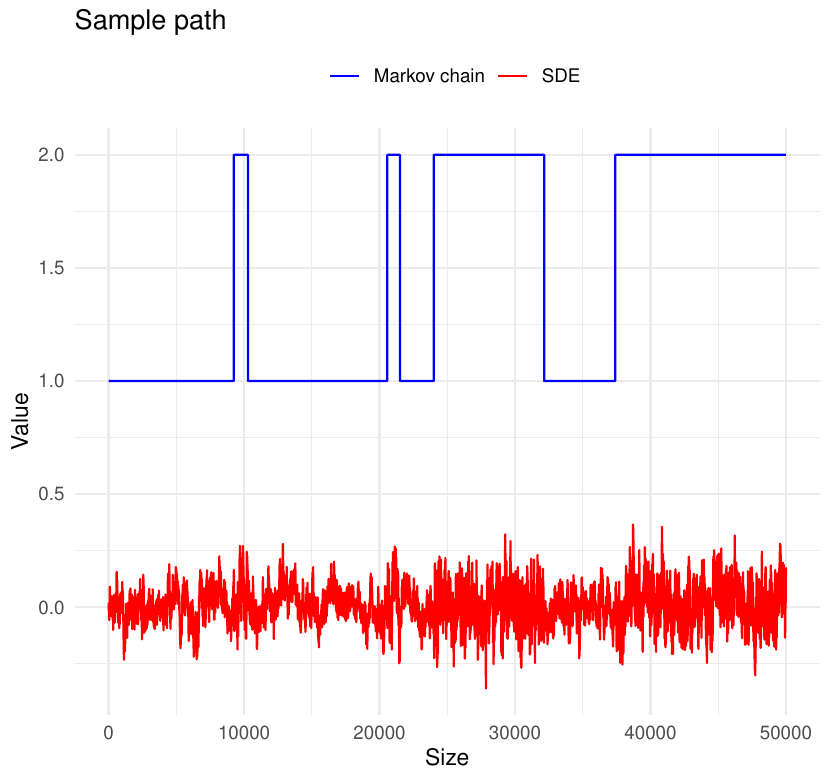}
	\caption{Sample paths of SDE solution $X$ and associated Markov chain $Z$.}
	\label{fig:samplepath_1dou}
\end{figure}

Here, we describe the method used to simulate high-frequency sample $\{(X_{t_j},Z_{t_j})\}_{j=0}^n$ with $t_j = jh$; recall that $h = T/n$ where $T$ denotes the terminal sampling time.
That is, generating one sample path is as follows.
\begin{itemize}
	\item First, we use the R package \texttt{spuRs} to generate Markov chain $Z$ (see \cite{spuRs_book} for details) with a smaller step size $\delta = h/10$. 
	\item Using the data $(Z_{\delta l})_{l=0}^{10 n}$, we then apply the Euler scheme to generate the data for $X$:
	\begin{equation*}
		X_{s_l} = -\beta_{Z_{s_{l-1}}}X_{s_{l-1}}\delta + \sigma_{Z_{s_{l-1}}}
		(w_{s_l}-w_{s_{l-1}}),
	\end{equation*}
	where $s_l = l\delta$ 
	for $l\in \{0,1,...,10n\}$; this internally generates $\{(X_{s_l},Z_{s_l})\}_{l=0}^{10n}$. 
	\item Next, we select a subsequence of the data $\{(X_{t_j},Z_{t_j})\}_{j=0}^{n}$ from $\{(X_{t_l},Z_{t_l})\}_{l=0}^{10n}$ with step size $h$, where $t_j = jh$. By thinning the sequence, we obtain high-frequency data $\{(X_{t_j},Z_{t_j})\}_{j=0}^{n}$.
\end{itemize}
We independently repeat the above procedure $M$ times.
Figure \ref{fig:samplepath_1dou} shows sample paths of $X$ with associated Markov chain $Z$ for $T=500$ and $n=50000$.

\begin{figure}[htp]
	\centering
	\includegraphics[width=\textwidth]{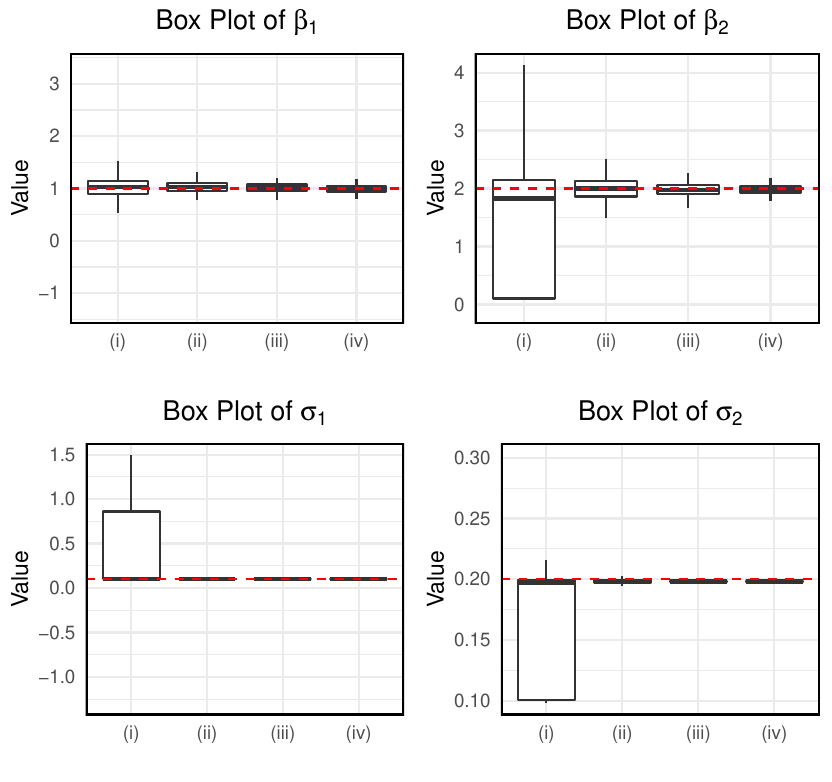}
	\caption{Boxplots of the estimators for $h=0.01$ with four schemes:(i) $T=100,n=10000$; (ii) $T=300,n=30000$;
		(iii) $T=500,n=50000$; (iv) $T=1000, n=100000$. The red dashed line indicates the true value of the parameters.}
	\label{fig:box1}
\end{figure}

\begin{figure}[htp]
	\centering
	\includegraphics[width=\textwidth]{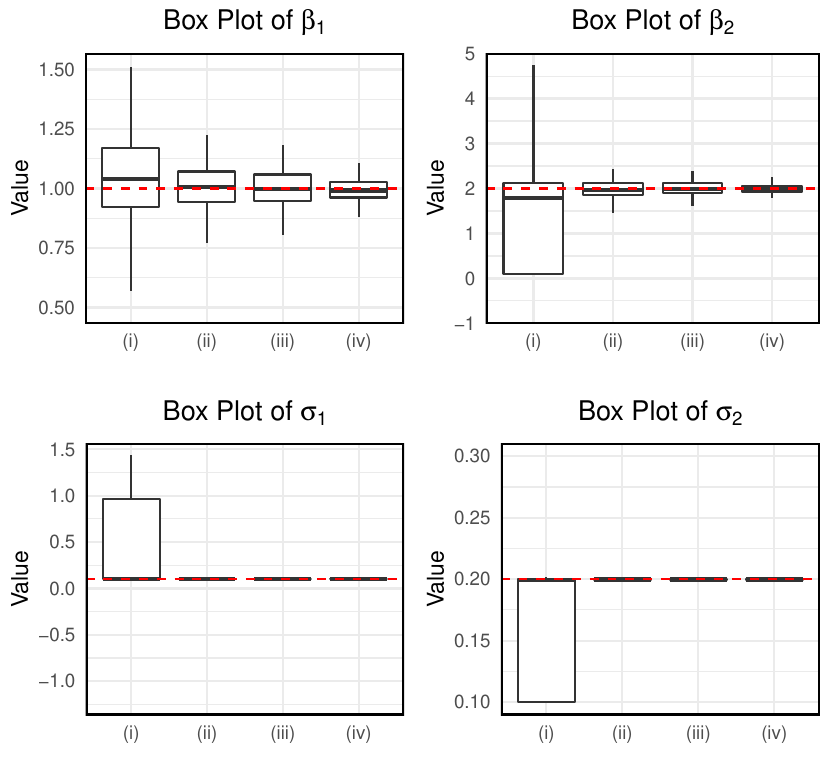}
	\caption{Boxplots of the estimators for $h=0.001$ with four schemes:(i) $T=100,n=100000$; (ii) $T=300,n=300000$;
		(iii) $T=500,n=500000$; (iv) $T=1000, n=1000000$. The red dashed line indicates the true value of the parameters.}
	\label{fig:box2}
\end{figure}

\begin{figure}[htbp]
	\centering
	\begin{minipage}{0.45\textwidth}
		\centering
		\includegraphics[width=\textwidth]{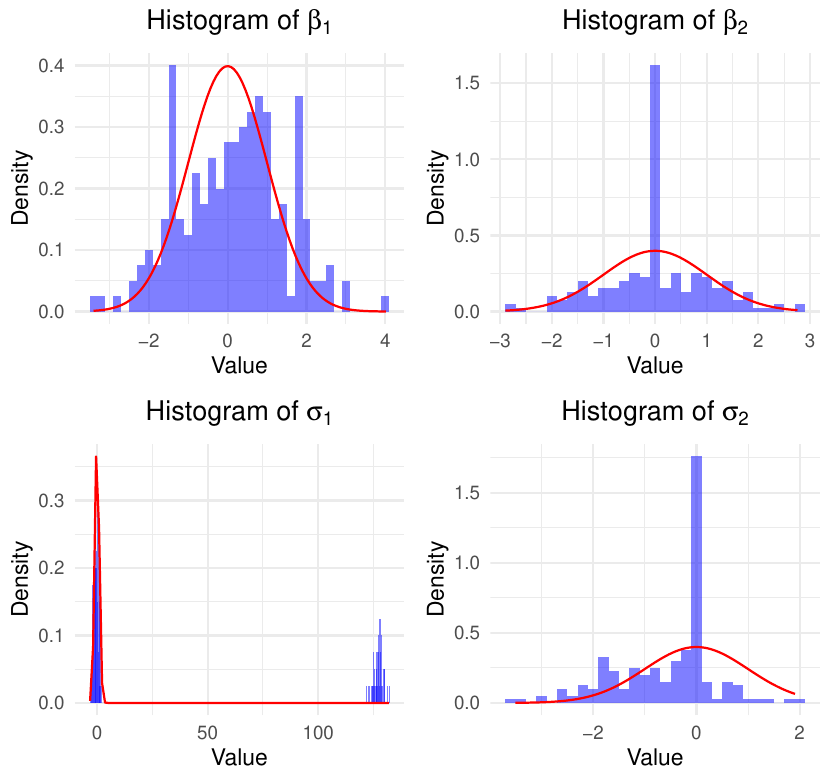}
		\subcaption{$T = 100$, $n = 10000$}
		\label{fig:dataset1}
	\end{minipage}
	\hfill
	\begin{minipage}{0.45\textwidth}
		\centering
		\includegraphics[width=\textwidth]{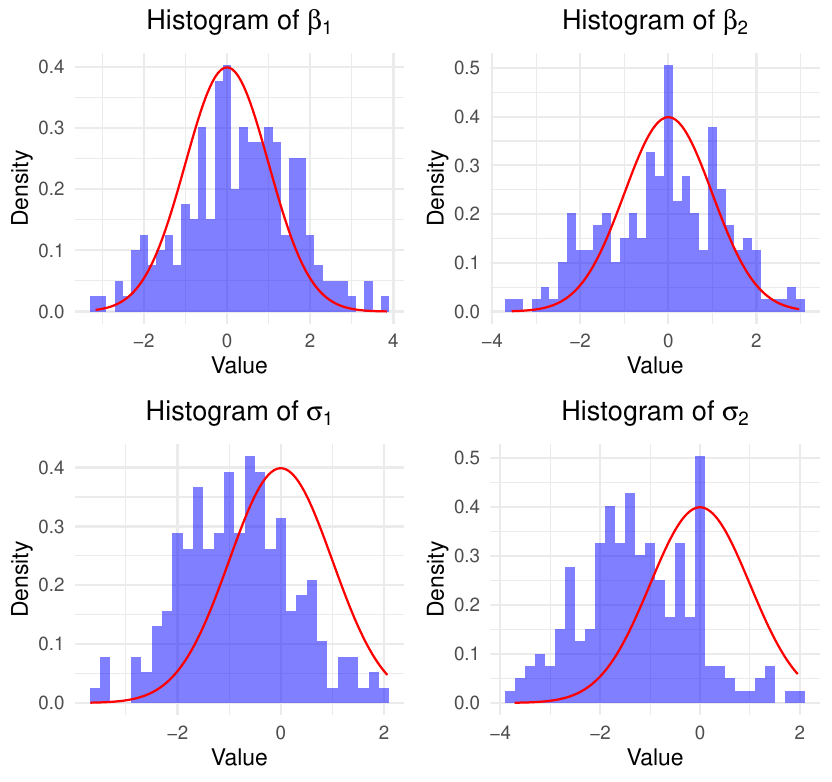}
		\subcaption{$T = 300$, $n = 30000$}
		\label{fig:dataset2}
	\end{minipage}
	
	\begin{minipage}{0.45\textwidth}
		\centering
		\includegraphics[width=\textwidth]{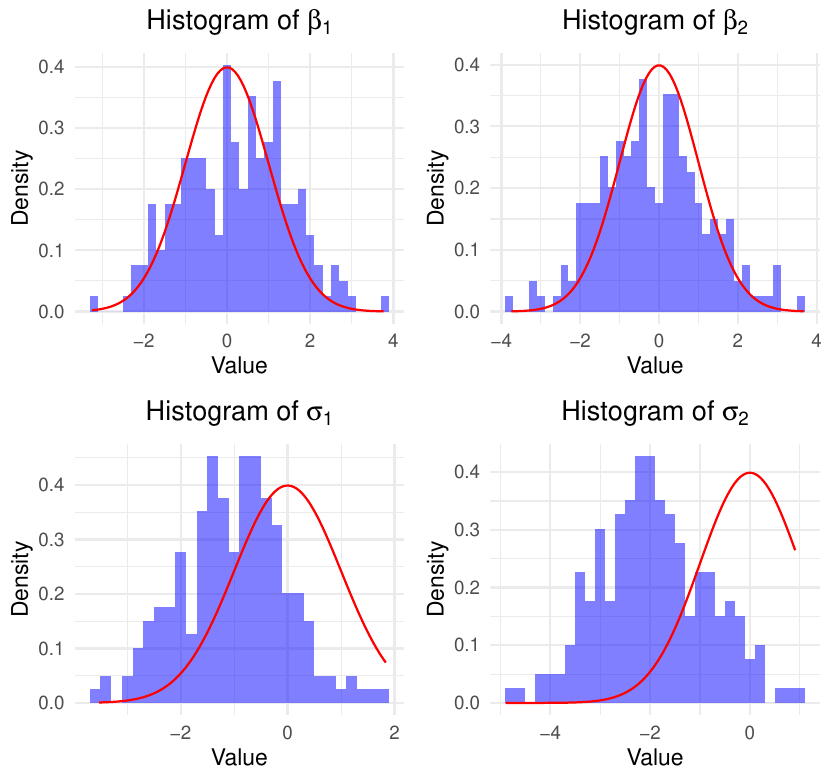}
		\subcaption{$T = 500$, $n = 50000$}
		\label{fig:dataset3}
	\end{minipage}
	\hfill
	\begin{minipage}{0.45\textwidth}
		\centering
		\includegraphics[width=\textwidth]{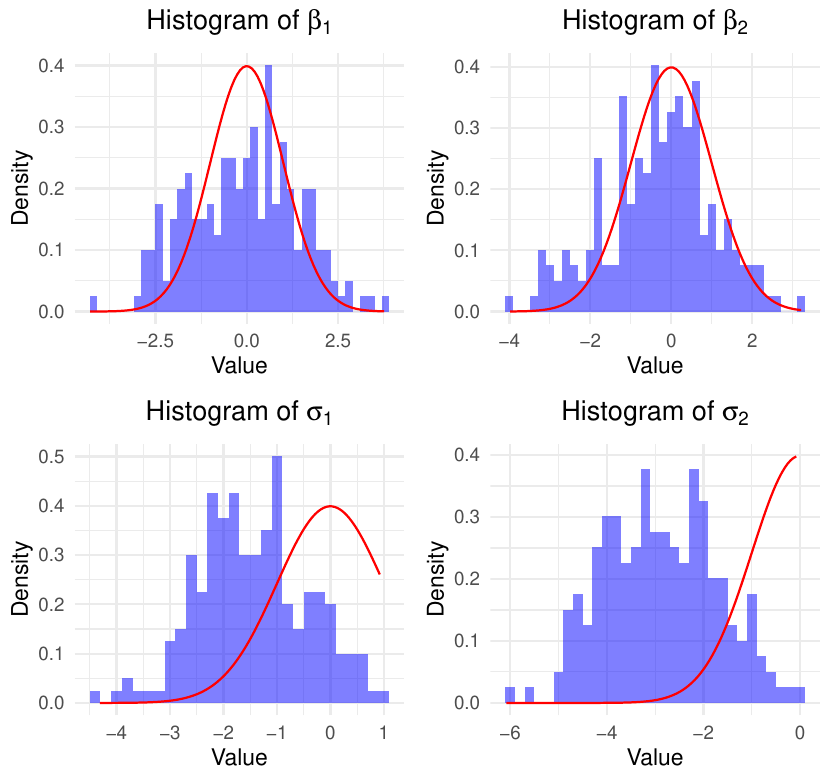}
		\subcaption{$T = 1000$, $n = 100000$}
		\label{fig:dataset4}
	\end{minipage}
	
	\caption{Histograms of the standardized estimators for $h = 0.01$. The red curve indicates the standard normal density.}
	\label{fig:figures1}
\end{figure}

\begin{figure}[htbp]
	\centering
	\begin{minipage}{0.45\textwidth}
		\centering
		\includegraphics[width=\textwidth]{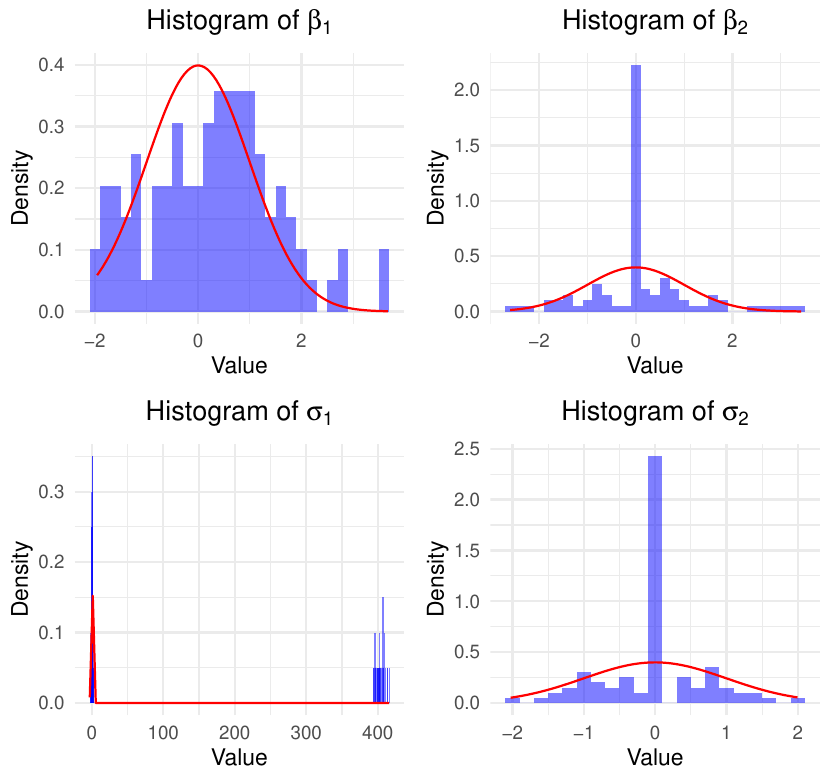}
		\subcaption{$T = 100$, $n = 100000$}
		\label{fig:dataset5}
	\end{minipage}
	\hfill
	\begin{minipage}{0.45\textwidth}
		\centering
		\includegraphics[width=\textwidth]{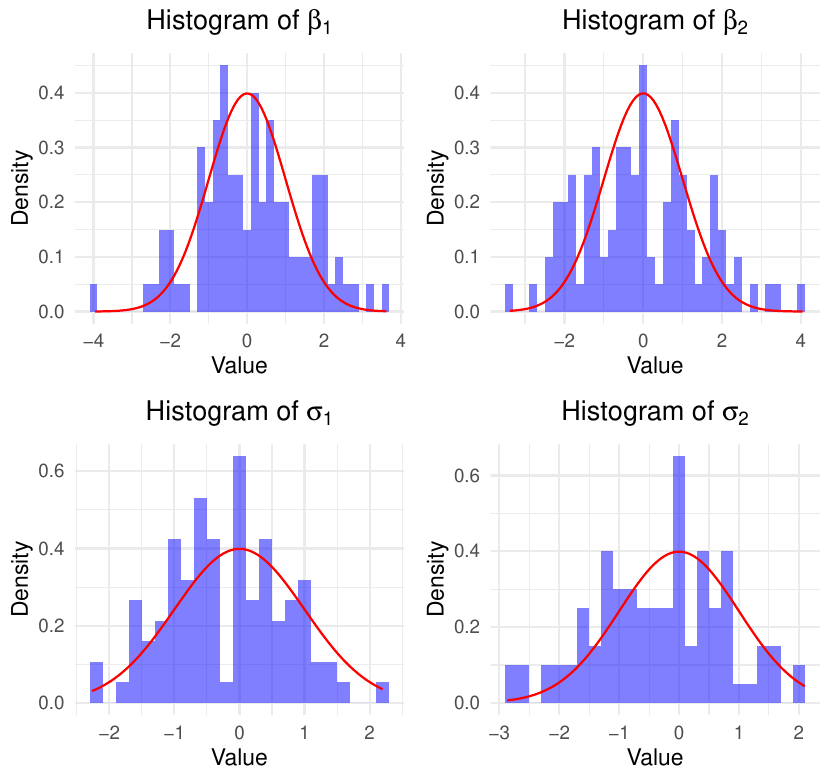}
		\subcaption{$T = 300$, $n = 300000$}
		\label{fig:dataset6}
	\end{minipage}
	
	\begin{minipage}{0.45\textwidth}
		\centering
		\includegraphics[width=\textwidth]{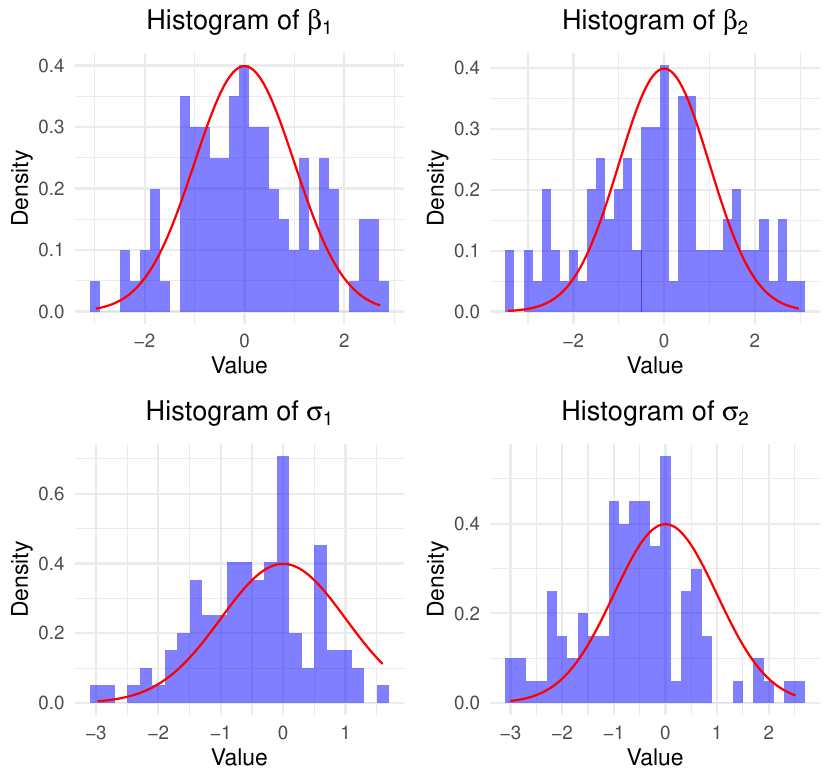}
		\subcaption{$T = 500$, $n = 500000$}
		\label{fig:dataset7}
	\end{minipage}
	\hfill
	\begin{minipage}{0.45\textwidth}
		\centering
		\includegraphics[width=\textwidth]{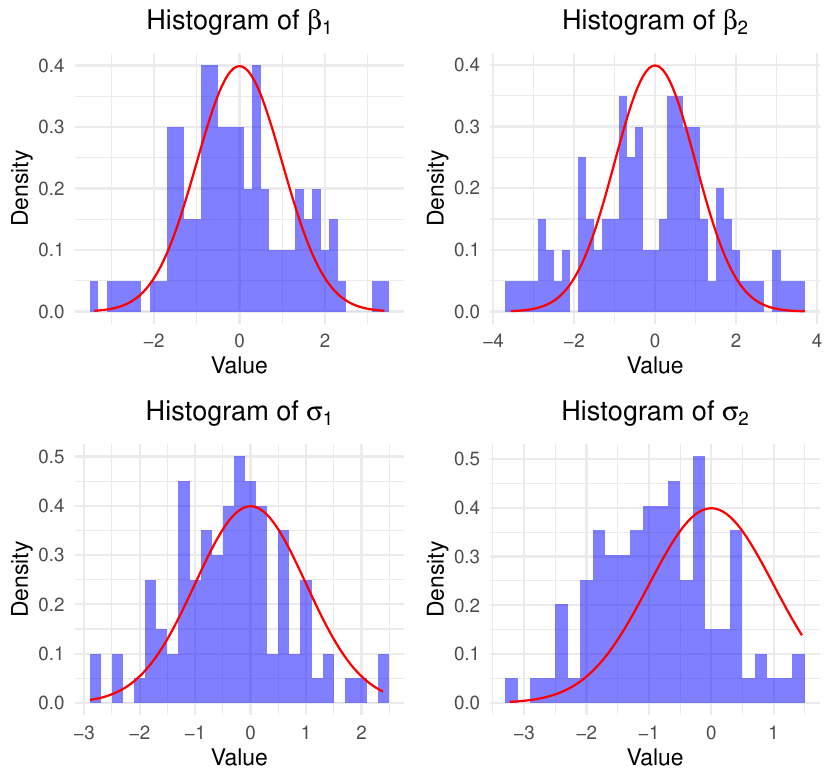}
		\subcaption{$T = 1000$, $n = 1000000$}
		\label{fig:dataset8}
	\end{minipage}
	
	\caption{Histograms of the standardized estimators for $h = 0.001$. The red curve indicates the standard normal density.}
	\label{fig:figures2}
\end{figure}

We set $M=200$ and $h=0.01$, and consider the following four schemes:
(i) $T=100,n=10000$; (ii) $T=300,n=30000$; (iii) $T=500,n=50000$; (iv) $T=1000, n=100000$. 
The estimators of $\theta$ and $Q$ are computed through formulas given in \eqref{estor:ou_1}, \eqref{estor:ou_2} and \eqref{eq:estQn}. Table \ref{tab:mean_sd_1} and Table \ref{tab:mean_sd_3} summarize the empirical mean and standard deviation (Std. Dev.) of these estimators. Figure \ref{fig:box1} shows the box plots of the estimators of $\theta$.
To illustrate the asymptotic normality stated in Theorem  \ref{thm:AN}, we present the standardized estimators for $\beta_i$ and $\sigma_i$ for each $i\in\{1,2\}$. 
By \eqref{hm:ex.eq-1} and \eqref{hm:ex.eq-2} in Section \ref{ex:ou}, 
they are given by
\begin{equation*}
	\sqrt{\frac{\hat{\pi}_i}{\hat{\beta}_i}}\sqrt{nh}(\hat{\beta}_i-\beta_i^\star),
	\qquad
	\sqrt{\frac{2\hat{\pi}_i}{(\hat{\sigma}_i)^2}}\sqrt{n}(\hat{\sigma}_i-\sigma_i^\star).
\end{equation*}
Figure \ref{fig:figures1} shows the histograms of the standardized estimators for each scheme.
Analogously, we generate $M=100$ sample paths for $(X_{t_j},Z_{t_j})$ for each of these four schemes, each with a step size of $h=0.001$: 
(i) $T=100,n=100000$; (ii) $T=300,n=300000$;
(iii) $T=500,n=500000$; (iv) $T=1000, n=1000000$. 
The estimators are computed as before. 
The mean and standard deviation of these computed estimators are presented in Table \ref{tab:mean_sd_2} and Table \ref{tab:mean_sd_4}. Figure \ref{fig:box2} shows the boxplots of the estimators for $\theta$. The corresponding histograms of the standardized estimators are shown in Figure \ref{fig:figures2}.

By examining the aforementioned simulation results, we observed the following:
\begin{itemize}
	\item Table \ref{tab:mean_sd_1} shows that the first scheme yields poorer estimates in comparison to the other three schemes. The second scheme generally offers good estimates, although it performs poorly for $\sigma_1$. The latter two schemes consistently produce accurate estimates.   In Table \ref{tab:mean_sd_2}, the latter two schemes continue to yield accurate estimates, whereas the first and second schemes provide comparatively less accurate estimates. 
	For both Table \ref{tab:mean_sd_1} and Table \ref{tab:mean_sd_2}, the estimated results improve as the terminal time $T$ and sample size $n$ increase.
	Looking at Table \ref{tab:mean_sd_3} and Table \ref{tab:mean_sd_4}, we observe that the estimates for $q_{21}$ and $q_{22}$ are poor when $T$ and $n$ are small, but improve as $T$ and $n$ increase.
	
	\item It is observed from the boxplots \ref{fig:box1} and \ref{fig:box2} that the estimates converge towards the true value in all plots as $T$ and $n$ increase. It is evident in scheme (i) where $T$ is relatively small, the range between the maximum and minimum estimates is large, indicating high variance. As $T$ increases, the estimates become more concentrated.  
	
	\item 
	By Theorem \ref{thm:AN}, histograms should fit the red standard normal curve well when $T$ is sufficiently large and $h$ is sufficiently small. 
	The performance of the histograms improves as the step size $h$ decreases.  This trend is evident in Figure \ref{fig:dataset3} compared to Figure \ref{fig:dataset7} and Figure \ref{fig:dataset4} compared to Figure \ref{fig:dataset8}. Overall, the estimates for $\beta$ show better performance than those for $\sigma$ under standardization. Specifically, for $\sigma$, the histograms exhibit a poor fit to the standard normal curve in Figures \ref{fig:dataset3} and \ref{fig:dataset4}. This may be attributed to the sensitivity of the condition $nh \to \infty$ required for asymptotic properties to hold, particularly regarding the parameter $\sigma$.
\end{itemize}

In simulations, although we theoretically assume the irreducibility of the Markov chain $Z_t$, meaning that it will range over all states when $T$ is sufficiently large, a small terminal time $T$ may not allow $Z_t$ to range over all states in $M$ sample paths. This can affect certain estimates in standardization, such as those of $\pi$, and ultimately lead to poor estimates and histogram fits. Therefore, as shown by our theoretical results, a sufficiently long terminal time and a sufficiently small step size are necessary to achieve reasonably accurate estimates.

\begin{table}[htbp]
	\centering
	\caption{The mean and the standard deviation (Std. Dev.) of the estimators with true values $\theta=(1,2,0.1,0.2)$ and $h=0.01$.}
	\begin{adjustbox}{max width=\textwidth}
		\begin{tabular}{ccccccccc}
			\toprule
			& \multicolumn{2}{c}{$T=100, n=10000$} & \multicolumn{2}{c}{$T=300, n=30000$} & \multicolumn{2}{c}{$T=500, n=50000$} & \multicolumn{2}{c}{$T=1000, n=100000$}  \\
			\cmidrule(lr){2-3} \cmidrule(lr){4-5} \cmidrule(lr){6-7} \cmidrule(lr){8-9}
			Estimator & Mean & Std. Dev. & Mean & Std. Dev. & Mean & Std. Dev. & Mean & Std. Dev. \\
			\midrule
			$\hat{\beta}_{1,n}$ & 1.130 & 1.018 & 1.027 & 0.125 & 1.018 & 0.088 & 0.999 & 0.072 \\
			$\hat{\beta}_{2,n}$ & 1.546 & 1.192 & 1.913 & 0.452 & 2.002 & 0.374 & 1.981 & 0.092 \\
			$\hat{\sigma}_{1,n}$ & 0.367 & 0.422 & 0.141 & 0.189 & 0.104 & 0.999 & 0.100 & 0.000 \\
			$\hat{\sigma}_{2,n}$ & 0.169 & 0.045 & 0.194 & 0.021 & 0.198 & 0.072 & 0.198 & 0.001 \\
			\bottomrule
		\end{tabular}
	\end{adjustbox}
	\label{tab:mean_sd_1}
\end{table}

\begin{table}[htbp]
	\centering
	\caption{The mean and the standard deviation (Std. Dev.) of the estimators with true values $\theta=(1,2,0.1,0.2)$ and $h=0.001$.}
	\begin{adjustbox}{max width=\textwidth}
		\begin{tabular}{ccccccccc}
			\toprule
			& \multicolumn{2}{c}{$T=100, n=100000$} & \multicolumn{2}{c}{$T=300, n=300000$} & \multicolumn{2}{c}{$T=500,n=500000$} & \multicolumn{2}{c}{$T=1000, n=1000000$} \\
			\cmidrule(lr){2-3} \cmidrule(lr){4-5} \cmidrule(lr){6-7} \cmidrule(lr){8-9} 
			Estimator & Mean & Std. Dev. & Mean & Std. Dev. & Mean & Std. Dev. & Mean & Std. Dev. \\
			\midrule
			$\hat{\beta}_{1,n}$ & 1.042 & 0.187 & 1.018 & 0.122 &  1.011 & 0.103 & 0.999 & 0.065 \\
			$\hat{\beta}_{2,n}$ & 1.372 & 1.157 & 1.879 & 0.526& 1.970 & 0.266 & 2.006 & 0.118 \\
			$\hat{\sigma}_{1,n}$ & 0.459 & 0.460 & 0.153 & 0.213& 0.110 & 0.102 & 0.100 & 0.000 \\
			$\hat{\sigma}_{2,n}$ & 0.161 & 0.049 & 0.194 & 0.024& 0.199 & 0.010 & 0.200 & 0.000 \\
			\bottomrule
		\end{tabular}
	\end{adjustbox}
	\label{tab:mean_sd_2}
\end{table}

\begin{table}[htbp]
	\centering
	\caption{The mean and the standard deviation (Std. Dev.) of the estimators with true values $Q =
		\begin{pmatrix}
			-0.01 & 0.01 \\
			0.01 & -0.01 \\
		\end{pmatrix}$ and $h=0.01$.}
	\begin{adjustbox}{max width=\textwidth}
		\begin{tabular}{ccccccccc}
			\toprule
			& \multicolumn{2}{c}{$T=100, n=10000$} & \multicolumn{2}{c}{$T=300, n=30000$} & \multicolumn{2}{c}{$T=500, n=50000$} & \multicolumn{2}{c}{$T=1000, n=100000$}  \\
			\cmidrule(lr){2-3} \cmidrule(lr){4-5} \cmidrule(lr){6-7} \cmidrule(lr){8-9}
			Estimator & Mean & Std. Dev. & Mean & Std. Dev. & Mean & Std. Dev. & Mean & Std. Dev. \\
			\midrule
			$\hat{q}_{11}^{(n)}$ & -0.0095 & 0.0172 & -0.0099 & 0.0097 & -0.0103 & 0.0071 & -0.0109 & 0.0057 \\
			$\hat{q}_{12}^{(n)}$ & 0.0095 & 0.0172 & 0.0099 & 0.0097 & 0.0103 & 0.0071 & 0.0109 & 0.0057 \\
			$\hat{q}_{21}^{(n)}$ & 0.0664 & 0.1463 & 0.0280 & 0.0645 & 0.0180 & 0.0387 & 0.0121 & 0.0061 \\
			$\hat{q}_{22}^{(n)}$ & -0.0664 & 0.1463 & -0.0280 & 0.0645 & -0.0180 & 0.0387 & -0.0121 & 0.0061 \\
			\bottomrule
		\end{tabular}
	\end{adjustbox}
	\label{tab:mean_sd_3}
\end{table}

\begin{table}[htbp]
	\centering
	\caption{The mean and the standard deviation (Std. Dev.) of the estimators with true values $Q =
		\begin{pmatrix}
			-0.01 & 0.01 \\
			0.01 & -0.01 \\
		\end{pmatrix}$ and $h=0.001$.}
	\begin{adjustbox}{max width=\textwidth}
		\begin{tabular}{ccccccccc}
			\toprule
			& \multicolumn{2}{c}{$T=100, n=100000$} & \multicolumn{2}{c}{$T=300, n=300000$} & \multicolumn{2}{c}{$T=500, n=500000$} & \multicolumn{2}{c}{$T=1000, n=1000000$}  \\
			\cmidrule(lr){2-3} \cmidrule(lr){4-5} \cmidrule(lr){6-7} \cmidrule(lr){8-9}
			Estimator & Mean & Std. Dev. & Mean & Std. Dev. & Mean & Std. Dev. & Mean & Std. Dev. \\
			\midrule
			$\hat{q}_{11}^{(n)}$ & -0.0133 & 0.0273 & -0.0098 & 0.0091 & -0.0107 & 0.0123 & -0.0102 & 0.0057 \\
			$\hat{q}_{12}^{(n)}$ & 0.0133 & 0.0273 & 0.0098 & 0.0091 & 0.0107 & 0.0123 & 0.0102 & 0.0057 \\
			$\hat{q}_{21}^{(n)}$ & 0.0643 & 0.1783 & 0.0270 & 0.0343 & 0.0163 & 0.0247 & 0.0126 & 0.0058 \\
			$\hat{q}_{22}^{(n)}$ & -0.0643 & 0.1783 & -0.0270 & 0.0343 & -0.0163 & 0.0247 & -0.0126 & 0.0058 \\
			\bottomrule
		\end{tabular}
	\end{adjustbox}
	\label{tab:mean_sd_4}
\end{table}

\section{Proofs}
\label{sec:proofs}

\subsection{Preliminary lemmas}

Let $C$ represent positive constants, the values of which may vary from one context to another. We use the notation $R_{j-1}(\theta)$ for 
\begin{equation*}
	\sup_\theta |R_{j-1}(\theta)|\le C(1+|X_{t_{j-1}}|)^C;
\end{equation*} we simply write $R_{j-1}$ if in particular $R_{j-1}$ does not depend on $\theta$. We also use the shorthand $\mathbb{E}_{\theta}^{j-1}(\cdot)= \mathbb{E}_{\theta}(\cdot|\mathcal{F}_{t_{j-1}})$.

From the ergodic theorem \eqref{eq:c-ergodic}, we can deduce the following discrete-time version:

\begin{lemma}
	Under Assumption \ref{ass:smoothness} and \ref{Ass:finite moment},  for $\rho \in \{\alpha, \gamma\}$,
	let $f(x,i,\rho_i)$ be a real-valued measurable function that is twice continuously differentiable with respect to the first variable and three times continuously differentiable with respect to the third variable, satisfying
	\begin{equation}
		\sup _{(x,\rho_i) \in \mathbb{R} \times \Theta_{\rho_i}} \frac{1}{1+|x|^C} \left|\partial_{\theta_i}^{k} \partial_{x}^{l} f(x, i, \rho_i)\right| < \infty
		\label{lem:con:smooth}
	\end{equation} where $k \in \{0,1,2,3\}$ and $l \in \{0,1,2\}$. Then 
	\begin{equation}
		\frac{1}{n} \sum_{j=1}^{n} f(X_{t_{j-1}},Z_{t_{j-1}},\rho_{Z_{t_{j-1}}}) \xrightarrow{p} \sum_{i =1}^{N}\int_{\mathbb{R}}f(x,i,\rho_i)\nu_{\theta^{\star}}(dx,i).
	\end{equation}
	\label{lem:disergodic}
\end{lemma}
\begin{proof}
	Note that
	\begin{align*}
		& \mathbb{E}_{\theta^{\star}}\left(\frac{1}{nh}\sum_{j=1}^{n} \int_{t_{j-1}}^{t_j}|f(X_s,Z_s,\rho_{Z_s})-f(X_{t_{j-1}},Z_{t_{j-1}},\rho_{Z_{t_{j-1}}})|ds\right) 
		\\
		&\leq \frac{1}{nh}\sum_{j=1}^{n} \int_{t_{j-1}}^{t_j}\mathbb{E}_{\theta^{\star}}\left(|f(X_s,Z_s,\rho_{Z_s})-f(X_{t_{j-1}},Z_s,\rho_{Z_s})|\right)ds 
		\\
		&\quad + \frac{1}{nh}\sum_{j=1}^{n} \int_{t_{j-1}}^{t_j}\mathbb{E}_{\theta^{\star}}\left(|f(X_{t_{j-1}},Z_s,\rho_{Z_s})-f(X_{t_{j-1}},Z_{t_{j-1}},\rho_{Z_{t_{j-1}}})|\right)ds
		\\
		&=: \mathcal{T}_1 +\mathcal{T}_2
	\end{align*}
	We first consider the term $\mathcal{T}_1$. It is important to note that  \cite[Theorem 3.23]{mao2006stochastic} implies that for $p \geq 2$ and $s \in [t_{j-1},t_j]$,
	\begin{equation*}
		\mathbb{E}_{\theta^{\star}}\left(|X_s-X_{t_{j-1}}|^p\right) = O(h^{\frac{p}{2}}).
	\end{equation*}
	Applying H\"{o}lder's inequality, condition \eqref{lem:con:smooth} and Assumption \ref{Ass:finite moment}, we have
	\begin{align*}
		& \frac{1}{nh}\sum_{j=1}^{n} \int_{t_{j-1}}^{t_j}\mathbb{E}_{\theta^{\star}}\left(|f(X_s,Z_s,\rho_{Z_s})-f(X_{t_{j-1}},Z_s,\rho_{Z_s})|\right)ds 
		\\
		&\leq \frac{1}{nh}\sum_{j=1}^{n} \int_{t_{j-1}}^{t_j}\left(\mathbb{E}_{\theta^{\star}}\left(X_s-X_{t_{j-1}}\right)^2\right)^{\frac{1}{2}}
		\\
		&\quad \left(
		\mathbb{E}_{\theta^{\star}}\left(\int_{0}^{1}\partial_x f(X_{t_{j-1}} + u(X_s-X_{t_{j-1}}),Z_s,\rho_{Z_s})du\right)^2\right)^{\frac{1}{2}}ds
		\\
		&=O(\sqrt{h}).
	\end{align*}
	For the term $\mathcal{T}_2$, equation  \eqref{CTMCgenerator} implies
	\begin{equation*}
		\mathbb{P}_{\theta^{\star}}(Z_s \neq i |Z_{t_{j-1}}=i) = \sum_{k\in S, k\neq i}(q_{ik}(s-t_{j-1}) +o(s-t_{j-1})).
	\end{equation*}
	This together with condition \eqref{lem:con:smooth} and Assumption \ref{Ass:finite moment} gives the following estimates:
	\begin{align*}
		&\frac{1}{nh}\sum_{j=1}^{n} \int_{t_{j-1}}^{t_j}\mathbb{E}_{\theta^{\star}}\left(|f(X_{t_{j-1}},Z_s,\rho_{Z_s})-f(X_{t_{j-1}},Z_{t_{j-1}},\rho_{Z_{t_{j-1}}})|\right)ds 
		\\
		&\leq \frac{1}{nh}\sum_{j=1}^{n} \int_{t_{j-1}}^{t_j}C\mathbb{E}_{\theta^{\star}}\left(I_{\{Z_s \neq Z_{t_{j-1}}\}}\right)ds
		\\
		&= C\frac{1}{nh}\sum_{j=1}^{n} \int_{t_{j-1}}^{t_j}\mathbb{E}_{\theta^{\star}}\left(\mathbb{E}_{\theta^{\star}}\left(I_{\{Z_s \neq Z_{t_{j-1}}\}}|Z_{t_{j-1}}\right)\right)ds
		\\
		&= C\frac{1}{nh}\sum_{j=1}^{n} \int_{t_{j-1}}^{t_j}\mathbb{E}_{\theta^{\star}}\left(\sum_{i\in S}I_{\{Z_{t_{j-1}} = i\}}\mathbb{E}_{\theta^{\star}}\left(I_{\{Z_s \neq i\}}|Z_{t_{j-1}}=i\right)\right)ds
		\\
		&=C\frac{1}{nh}\sum_{j=1}^{n} \int_{t_{j-1}}^{t_j}\mathbb{E}_{\theta^{\star}}\left(\sum_{i\in S}I_{\{Z_{t_{j-1}} = i\}}\mathbb{P}_{\theta^{\star}}(Z_s \neq i |Z_{t_{j-1}}=i)\right)ds
		\\
		&= C\frac{1}{nh}\sum_{j=1}^{n} \int_{t_{j-1}}^{t_j}\mathbb{E}_{\theta^{\star}}\left(
		\sum_{i\in S}I_{\{Z_{t_{j-1}} = i\}}\sum_{k\in S, k\neq i}(q_{ik}(s-t_{j-1}) +o(s-t_{j-1}))
		\right)ds
		\\
		&= O(h).
	\end{align*}
	Combining the estimates on $\mathcal{T}_1$ and $\mathcal{T}_2$ with the property \eqref{eq:c-ergodic},  we obtain the desired result.  
\end{proof}

We present the following moment estimates for future use.
\begin{lemma}
	\label{lem:estimates}
	Under Assumption \ref{ass:smoothness}-\ref{Ass:finite moment}, we have
	\begin{align*}
		&\mathbb{E}_{\theta^{\star}}^{j-1}\left[(X_{t_j}-\mu_{j-1}(\alpha))^2\right] = \sigma^2(X_{t_{j-1}},Z_{t_{j-1}},\gamma_{Z_{t_{j-1}}}^{\star})h + h^{3/2} R_{j-1}(\theta),
		\\
		&\mathbb{E}_{\theta^{\star}}^{j-1}\left[X_{t_j}-\mu_{j-1}(\alpha^{\star})\right] = h^{3/2}R_{j-1},
		\\
		&\mathbb{E}_{\theta^{\star}}^{j-1}\left[(X_{t_j}-\mu_{j-1}(\alpha^{\star}))^3\right] = h^{5/2}R_{j-1},
		\\
		&\mathbb{E}_{\theta^{\star}}^{j-1}\left[(X_{t_j}-\mu_{j-1}(\alpha^{\star}))^4\right] = 3\sigma^4(X_{t_{j-1}},Z_{t_{j-1}},\gamma_{Z_{t_{j-1}}}^{\star})h^2 + h^3 R_{j-1}.
	\end{align*}
\end{lemma}

\begin{proof}
	We will only present the estimate for $\mathbb{E}_{\theta^{\star}}^{j-1}(X_{t_j}-\mu_{j-1}(\alpha))^2$. Note that within the expectation $\mathbb{E}_{\theta^{\star}}^{j-1}(X_{t_j}-\mu_{j-1}(\alpha^{\star}))$ and $\mathbb{E}_{\theta^{\star}}^{j-1}(X_{t_j}-\mu_{j-1}(\alpha^{\star}))^4$, the true parameter $\alpha^{\star}$ is used, which differs from the parameter in the first expectation. Keeping this distinction in mind, the estimate for $\mathbb{E}_{\theta^{\star}}^{j-1}(X_{t_j}-\mu_{j-1}(\alpha^{\star}))$ and $\mathbb{E}_{\theta^{\star}}^{j-1}(X_{t_j}-\mu_{j-1}(\alpha^{\star}))^4$ can be derived using a similar approach employed for the estimate of $\mathbb{E}_{\theta^{\star}}^{j-1}(X_{t_j}-\mu_{j-1}(\alpha))^2$. 
	We omit the related details. 
	
	Observe that
	\begin{align*}
		\mathbb{E}_{\theta^{\star}}^{j-1}\left[(X_{t_j}-\mu_{j-1}(\alpha))^2\right] =&
		\mathbb{E}_{\theta^{\star}}^{j-1}\Bigg[\Bigg(\int_{t_{j-1}}^{t_j}(b(X_s,Z_s,\alpha_{Z_s}^{\star})-b(X_{t_{j-1}},Z_{t_{j-1}},\alpha_{Z_{t_{j-1}}}))ds
		\\
		&+ \int_{t_{j-1}}^{t_j}(\sigma(X_s,Z_s,\gamma_{Z_s}^{\star})-\sigma(X_{t_{j-1}},Z_{t_{j-1}},\gamma_{Z_{t_{j-1}}}^{\star}))dw_s
		\\
		&+ \int_{t_{j-1}}^{t_j}\sigma(X_{t_{j-1}},Z_{t_{j-1}},\gamma_{Z_{t_{j-1}}}^{\star})dw_s\Bigg)^2\Bigg].
	\end{align*}
	We estimate all squared terms and mixed terms from the right-hand side of the above equation.
	By applying Assumption \ref{ass:smoothness}, \cite[Theorem 3.23]{mao2006stochastic}, and a similar argument as in the proof of Lemma \ref{lem:disergodic}, we observe that
	\begin{align*}
		\mathbb{E}_{\theta^{\star}}^{j-1}&\left[\left(b(X_s,Z_s,\alpha_{Z_s}^{\star})-b(X_{t_{j-1}},Z_{t_{j-1}},\alpha_{Z_{t_{j-1}}})\right)^2\right]
		\\
		\leq& 
		2\mathbb{E}_{\theta^{\star}}^{j-1}\left[\left(b(X_s,Z_s,\alpha_{Z_s}^{\star})-b(X_{t_{j-1}},Z_s,\alpha_{Z_s}^{\star})\right)^2 \right]
		\\
		&+ 2\mathbb{E}_{\theta^{\star}}^{j-1}\left[\left(b(X_{t_{j-1}},Z_s,\alpha_{Z_s}^{\star})-b(X_{t_{j-1}},Z_{t_{j-1}},\alpha_{Z_{t_{j-1}}})\right)^2\right]
		\\
		\leq&  h R_{j-1} + C\mathbb{E}_{\theta^{\star}}^{j-1}\left[\left(b(X_{t_{j-1}},Z_s,\alpha_{Z_s}^{\star})-b(X_{t_{j-1}},Z_s,\alpha_{Z_s})\right)^2\right]
		\\
		&+ C \mathbb{E}_{\theta^{\star}}^{j-1}\left[\left(b(X_{t_{j-1}},Z_s,\alpha_{Z_s})-b(X_{t_{j-1}},Z_{t_{j-1}},\alpha_{Z_{t_{j-1}}})\right)^2\right]
		\\
		\leq& h R_{j-1} + R_{j-1}(\theta) + R_{j-1}(\theta)\mathbb{E}_{\theta^{\star}}^{j-1}\left(I_{\{Z_s \neq Z_{t_{j-1}}\}}\right) 
		\\
		\leq& (1+h) R_{j-1}(\theta).
	\end{align*}
	Therefore, Jensen's inequality implies
	\begin{align*}
		&\mathbb{E}_{\theta^{\star}}^{j-1}\left(\int_{t_{j-1}}^{t_j}(b(X_s,Z_s,\alpha_{Z_s}^{\star})-b(X_{t_{j-1}},Z_{t_{j-1}},\alpha_{Z_{t_{j-1}}}))ds \right)^2
		\\
		&\leq h\int_{t_{j-1}}^{t_j}\mathbb{E}_{\theta^{\star}}^{j-1}\left[(b(X_s,Z_s,\alpha_{Z_s}^{\star})-b(X_{t_{j-1}},Z_{t_{j-1}},\alpha_{Z_{t_{j-1}}}))^2\right]ds 
		\\
		&= h^2 R_{j-1}(\theta).
	\end{align*}
	By applying similar procedures, it is easy to obtain that 
	\begin{align*}
		\mathbb{E}_{\theta^{\star}}^{j-1}\left[\left(\sigma(X_s,Z_s,\gamma_{Z_s}^{\star})-\sigma(X_{t_{j-1}},Z_{t_{j-1}},\gamma_{Z_{t_{j-1}}}^{\star})\right)^2\right]
		= h R_{j-1},
	\end{align*}
	which implies
	\begin{align*}
		&\mathbb{E}_{\theta^{\star}}^{j-1}\left[\left(\int_{t_{j-1}}^{t_j}(\sigma(X_s,Z_s,\gamma_{Z_s}^{\star})-\sigma(X_{t_{j-1}},Z_{t_{j-1}},\gamma_{Z_{t_{j-1}}}^{\star}))dw_s\right)^2\right]
		\\
		&=
		\int_{t_{j-1}}^{t_j}\mathbb{E}_{\theta^{\star}}^{j-1}\left[\left(\sigma(X_s,Z_s,\gamma_{Z_s}^{\star})-\sigma(X_{t_{j-1}},Z_{t_{j-1}},\gamma_{Z_{t_{j-1}}}^{\star})\right)^2 \right]ds
		\\
		&= h^2 R_{j-1}.
	\end{align*}
	Obviously,
	\begin{equation*}
		\mathbb{E}_{\theta^{\star}}^{j-1}\left[\left(\int_{t_{j-1}}^{t_j}\sigma(X_{t_{j-1}},Z_{t_{j-1}},\gamma_{Z_{t_{j-1}}}^{\star})dw_s\right)^2\right] = \sigma^2(X_{t_{j-1}},Z_{t_{j-1}},\gamma_{Z_{t_{j-1}}}^{\star})h.
	\end{equation*}
	It is easy to estimate all other mixed terms by using the H\"{o}lder inequality and the above observations. We list the estimates of mixed terms below.
	\begin{align*}
		&\mathbb{E}_{\theta^{\star}}^{j-1} \left( \int_{t_{j-1}}^{t_j} A_s ds \int_{t_{j-1}}^{t_j} B_s dw_s \right) = h^2 R_{j-1}(\theta), \\
		&\mathbb{E}_{\theta^{\star}}^{j-1} \left( \int_{t_{j-1}}^{t_j} A_s ds \int_{t_{j-1}}^{t_j} \sigma(X_{t_{j-1}},Z_{t_{j-1}},\gamma_{Z_{t_{j-1}}}^{\star}) dw_s \right) = h^{3/2} R_{j-1}(\theta), \\
		&\mathbb{E}_{\theta^{\star}}^{j-1} \left( \int_{t_{j-1}}^{t_j} B_s ds \int_{t_{j-1}}^{t_j} \sigma(X_{t_{j-1}},Z_{t_{j-1}},\gamma_{Z_{t_{j-1}}}^{\star}) dw_s \right) = h^{3/2} R_{j-1}.
	\end{align*}
	where
	\begin{align*}
		&A_s = b(X_s, Z_s, \alpha_{Z_s}^{\star}) - b(X_{t_{j-1}}, Z_{t_{j-1}}, \alpha_{Z_{t_{j-1}}}), \\
		&B_s = \sigma(X_s, Z_s, \gamma_{Z_s}^{\star}) - \sigma(X_{t_{j-1}}, Z_{t_{j-1}}, \gamma_{Z_{t_{j-1}}}^{\star}).
	\end{align*}
	Combining the above observations, we obtain the desired result
	\begin{align*}
		\mathbb{E}_{\theta^{\star}}^{j-1}\left[(X_{t_j}-\mu_{j-1}(\alpha))^2\right] = \sigma^2(X_{t_{j-1}},Z_{t_{j-1}},\gamma_{Z_{t_{j-1}}}^{\star})h + h^{3/2} R_{j-1}(\theta).
	\end{align*}
\end{proof}

\subsection{Proof of Theorem \ref{thm:AN}}

\subsubsection{Consistency}

The consistency of $\hat{\theta}_n$ can be established by employing a method similar to that presented in \cite{Kes97}, leveraging the argmax theorem twice. See \cite{van2000} for details. 

First, we consider the parameter $\gamma$. Define
\begin{align*}
	&\mathbb{F}_{1,n}(\gamma):= \frac{1}{n}\left(\mathbb{H}_n(\alpha,\gamma)-\mathbb{H}_n(\alpha,\gamma^{\star})\right),
	\\
	&\mathbb{F}_{1}(\gamma):= -\frac{1}{2}\sum_{i =1}^{N}\int_{\mathbb{R}} \left(\log\left(\frac{\sigma^2(x,i,\gamma_i)}{\sigma^2(x,i,\gamma_i^{\star})}\right) + \frac{\sigma^2(x,i,\gamma_i^{\star})}{\sigma^2(x,i,\gamma_i)}-1 \right)\nu_{\theta^{\star}}(dx,i).
\end{align*}
Observe that 
\begin{equation*}
	\log\left(\frac{\sigma^2(x,i,\gamma_i)}{\sigma^2(x,i,\gamma_i^{\star})}\right) + \frac{\sigma^2(x,i,\gamma_i^{\star})}{\sigma^2(x,i,\gamma_i)}-1  \geq 0.
\end{equation*}
The function $\mathbb{F}_{1}(\gamma)$ is maximized only when $\sigma^2(x,i,\gamma_i)=\sigma^2(x,i,\gamma_i^{\star})$ for almost sure $x$ and each $i \in S$, therefore Assumption \ref{ass:identi} yields that $\mathop{\rm argmax}_{\gamma}\mathbb{F}_{1} = \{\gamma^{\star}\}$. 
Then, the argmax theorem concludes the consistency $\hat{\gamma}_n \xrightarrow{p} \gamma^{\star}$ if we show the uniform convergence in probability:
\begin{equation*}
	\sup_{\gamma} |\mathbb{F}_{1,n}(\gamma)-\mathbb{F}_{1}(\gamma)| \xrightarrow{p} 0.
\end{equation*}
To establish the desired result, we begin by examining the pointwise convergence. From By \ref{lem:disergodic}, Lemma \ref{lem:estimates}, and Assumption \ref{Ass:finite moment}, we have
\begin{align*}
	&\mathbb{F}_{1,n}(\gamma)
	\\
	 &= \frac{1}{n} \sum_{j=1}^{n} \frac{1}{2} \Bigg\{ \log \left(\frac{\sigma_{j-1}^2(\gamma^{\star})}{\sigma_{j-1}^2(\gamma)}\right) + \frac{1}{h}(X_{t_j}-\mu_{j-1}(\alpha))^2\left(\frac{1}{\sigma_{j-1}^2(\gamma^{\star})}-\frac{1}{\sigma_{j-1}^2(\gamma)}\right) \Bigg\}
	\\
	&= \frac{1}{n} \sum_{j=1}^{n} \frac{1}{2} \log \left(\frac{\sigma_{j-1}^2(\gamma^{\star})}{\sigma_{j-1}^2(\gamma)}\right) + \frac{1}{n} \sum_{j=1}^{n} \frac{1}{2} \left(1-\frac{\sigma_{j-1}^2(\gamma^{\star})}{\sigma_{j-1}^2(\gamma)}\right) + O_p(\sqrt{h})
	\\
	&= \mathbb{F}_{1}(\gamma) + o_p(1).
\end{align*}
In the above derivation, the transition to the second line is justified by the following procedure.
Note that by Lemma \ref{lem:estimates},
\begin{align*}
	&\frac{1}{n}\sum_{j=1}^{n} \mathbb{E}_{\theta^{\star}}^{j-1}\left(\frac{(X_{t_j}-\mu_{j-1}(\alpha))^2}{h\sigma_{j-1}^2(\gamma^{\star})}\right) = 1 + \frac{1}{n}\sum_{j=1}^{n} h^{1/2}R_{j-1}(\theta^{\star}) = 1 + O_p(\sqrt{h}),
	\\
	&\frac{1}{n^2}\sum_{j=1}^{n} \mathbb{E}_{\theta^{\star}}^{j-1}\left(\frac{(X_{t_j}-\mu_{j-1}(\alpha))^4}{h^2\sigma_{j-1}^4(\gamma^{\star})}\right) = \frac{1}{n^2}\sum_{j=1}^{n} h R_{j-1}(\theta^{\star}) = o_p(1).
\end{align*}
Then, by using \cite[Lemma 9]{GenonCatalot1993OnTE} or \cite[Lemma 3.4]{kessler2012statistical}, we have
\begin{align*}
	\frac{1}{n}\sum_{j=1}^{n} \frac{(X_{t_j}-\mu_{j-1}(\alpha))^2}{h\sigma_{j-1}^2(\gamma^{\star})} = 1 + O_p(\sqrt{h}).
\end{align*} Similarly, by applying analogous reasoning, we deduce that
\begin{align*}
	\frac{1}{n}\sum_{j=1}^{n} \frac{(X_{t_j}-\mu_{j-1}(\alpha))^2}{h\sigma_{j-1}^2(\gamma)} = \frac{1}{n}\sum_{j=1}^{n} \frac{\sigma_{j-1}^2(\gamma^{\star})}{\sigma_{j-1}^2(\gamma)} + O_p(\sqrt{h}).
\end{align*}
This procedure will be used frequently in the following arguments without further mention. It follows that $\mathbb{F}_{1,n}(\gamma) \xrightarrow{p} \mathbb{F}_{1}(\gamma)$ for each $\gamma$. 

It remains to show the tightness of the sequence $\{\mathbb{F}_{1,n}(\gamma)\}$ as random functions taking values in $C(\Theta_{\gamma})$. Since $\mathbb{F}_{1,n}(\gamma)= O_p(1)$ for every $\gamma$ by H\"{o}lder type inequalities, the tightness can be established by checking the Kolmogorov tightness condition \cite{kunita1986tightness}: 
it suffices to prove that
there exist $\delta > 0$, $\beta > 0$, and $C>0$ such that for every $\gamma,\gamma'$,
\begin{equation}
	\sup_n \mathbb{E}_{\theta^{\star}}\left(\left|\mathbb{F}_{1,n}(\gamma)-\mathbb{F}_{1,n}(\gamma')\right|^{\delta}\right) \leq C \left|\gamma-\gamma'\right|^{N+\beta}.
	\label{eq:kol_tight_con}
\end{equation}
Since we have
\begin{equation*}
	\left|\mathbb{F}_{1,n}(\gamma)-\mathbb{F}_{1,n}(\gamma')\right|^{\delta} \leq \sup_{\gamma} \left|\partial_{\gamma}\mathbb{F}_{1,n}(\gamma)\right|^{\delta} \left|\gamma-\gamma'\right|^{\delta},
\end{equation*}
we take $\delta = N+\beta$ for $\beta>0$, entailing that it is sufficient to check
\begin{align*}
	\sup_n \mathbb{E}_{\theta^{\star}}\left(\sup_{\gamma} \left|\partial_{\gamma}\mathbb{F}_{1,n}(\gamma)\right|^{N+\beta}\right) < \infty.
\end{align*}
This immediately follows from the smoothness condition in Assumption \ref{ass:smoothness} and the finite moment condition in Assumption \ref{Ass:finite moment} by applying H\"{o}lder type inequalities. 
We have established the tightness of sequence $\{\mathbb{F}_{1,n}(\gamma)\}$, hence the uniform convergence $\sup_{\gamma} |\mathbb{F}_{1,n}(\gamma)-\mathbb{F}_{1}(\gamma)| \xrightarrow{p} 0$, followed by 
$\hat{\gamma}_n \xrightarrow{p} \gamma^{\star}$.

Next, we consider the consistency of $\hat{\alpha}_n$. Let $b_{j-1}(\alpha):=b(X_{t_{j-1}},Z_{t_{j-1}},\alpha_{Z_{t_{j-1}}})$.
We define 
\begin{align*}
	&\mathbb{F}_{2,n}(\alpha):= \frac{1}{nh}\left(\mathbb{H}_n(\alpha,\gamma)-\mathbb{H}_n(\alpha^{\star},\gamma)\right),
	\\
	&\mathbb{F}_{2}(\alpha):= -\sum_{i =1}^{N}\int_{\mathbb{R}} \frac{\left(b(x,i,\alpha_i^{\star})-b(x,i,\alpha_i)\right)^2}{2\sigma^2(x,i,\gamma_i)}\nu_{\theta^{\star}}(dx,i).
\end{align*}
Also, we have $\mathop{\rm argmax}_{\alpha}\mathbb{F}_{2}=\{\alpha^{\star}\}$.
Applying Lemma \ref{lem:disergodic}, Lemma \ref{lem:estimates} and Assumption \ref{Ass:finite moment}, we have
\begin{align*}
	&\mathbb{F}_{2,n}(\alpha) 
	\\
	&= \frac{1}{nh} \sum_{j=1}^{n} \frac{1}{2h \sigma_{j-1}^2(\gamma)}\left(\left(X_{t_j}-\mu_{j-1}(\alpha^{\star})\right)^2-\left(X_{t_j}-\mu_{j-1}(\alpha)\right)^2\right)
	\\
	&= \frac{1}{nh} \sum_{j=1}^{n} \Bigg\{ \frac{(X_{t_j}-\mu_{j-1}(\alpha^{\star}))(\mu_{j-1}(\alpha)-\mu_{j-1}(\alpha^{\star}))}{h\sigma_{j-1}^2(\gamma)}  - \frac{\left(\mu_{j-1}(\alpha^{\star})-\mu_{j-1}(\alpha)\right)^2}{2h\sigma_{j-1}^2(\gamma)} \Bigg\}
	\\
	&= -\frac{1}{n} \sum_{j=1}^{n} \frac{\left(b_{j-1}(\alpha^{\star})-b_{j-1}(\alpha)\right)^2}{2\sigma_{j-1}^2(\gamma)} + O_p(\sqrt{h})
	\\
	&= \mathbb{F}_{2}(\alpha) + o_p(1).
\end{align*}
In a similar approach as for $\{\mathbb{F}_{1,n}(\alpha)\}$, we show the tightness of sequence $\{\mathbb{F}_{2,n}(\alpha)\}$ as random functions in $C(\Theta_{\alpha})$. The H\"{o}dler's inequalities ensure $\mathbb{F}_{2,n}(\alpha)= O_p(1)$ for every $\alpha$. To verify the Kolmogorov tightness criterion for $\{\mathbb{F}_{1,n}(\alpha)\}$, we apply Lemma \ref{lem:estimates} and Assumption \ref{ass:smoothness} to conlude
\begin{align*}
	&\sup_n \mathbb{E}_{\theta^{\star}}\left(\left|\mathbb{F}_{2,n}(\alpha)-\mathbb{F}_{2,n}(\alpha')\right|^{2}\right)
	\\
	&= \sup_n \mathbb{E}_{\theta^{\star}}\left(\left|\frac{1}{nh} \sum_{j=1}^{n} \frac{b_{j-1}(\alpha)-b_{j-1}(\alpha')} {\sigma_{j-1}^2(\gamma)} \Bigg\{ (X_{t_j}-\mu_{j-1}(\alpha')) 
	- \frac{\mu_{j-1}(\alpha)-\mu_{j-1}(\alpha')}{2} \Bigg\}\right|^{2}\right)
	\\
	&\leq C \sup_n  \frac{1}{n^2h^2} \sum_{j=1}^{n} \mathbb{E}_{\theta^{\star}}\left( \left|\frac{(b_{j-1}(\alpha)-b_{j-1}(\alpha'))^2} {\sigma_{j-1}^4(\gamma)}
	\Bigg\{ (X_{t_j}-\mu_{j-1}(\alpha')) 
	- \frac{\mu_{j-1}(\alpha)-\mu_{j-1}(\alpha')}{2} \Bigg\}\right|^2
	\right)
	\\
	&\leq C \sup_n  \frac{1}{n^2h^2} \sum_{j=1}^{n} \mathbb{E}_{\theta^{\star}}\left( \frac{(b_{j-1}(\alpha)-b_{j-1}(\alpha'))^2} {\sigma_{j-1}^4(\gamma)}
	\Bigg\{ \mathbb{E}_{\theta^{\star}}^{j-1}(X_{t_j}-\mu_{j-1}(\alpha'))^2 
	+ \frac{(b_{j-1}(\alpha)-b_{j-1}(\alpha'))^2h^2}{4} \Bigg\}
	\right)
	\\
	&\leq C \sup_n \left(\frac{1}{n^2}+\frac{1}{n^2h}\right) \sum_{j=1}^{n} \mathbb{E}_{\theta^{\star}}\left(
	\left(1+|X_{t_{j-1}}|\right)^C \sup_{\alpha}|\partial_{\alpha}b_{j-1}(\alpha)|^2
	\right)|\alpha-\alpha'|^2
	\\
	&\leq C |\alpha-\alpha'|^2.
\end{align*}
This leads to the tightness of the sequence $\{\mathbb{F}_{2,n}(\alpha)\}$. 
It follows from $\mathbb{F}_{2,n}(\alpha) = \mathbb{F}_{2}(\alpha) + o_p(1)$ for each $\alpha$ and the tightness that $\sup_{\alpha}|\mathbb{F}_{2,n}(\alpha)-\mathbb{F}_{2}(\alpha)|\xrightarrow{p} 0$. The argmax theorem now concludes that $\hat{\alpha}_n \xrightarrow{p} \alpha^{\star}$, which combined with $\hat{\gamma}_n \xrightarrow{p} \gamma^{\star}$ gives $\hat{\theta}_n \xrightarrow{p} \theta^{\star}$.

\subsubsection{Asymptotic normality}

It follows by a Taylor expansion of $\partial_{\theta}\mathbb{H}_n(\hat{\theta}_n)$ around $\theta^{\star}$ that
\begin{align*}
	0 = \partial_{\theta}\mathbb{H}_n(\hat{\theta}_n) = D_n^{-1}\partial_{\theta}\mathbb{H}_n(\theta^{\star}) + \left(\int_0^1 D_n^{-1}\partial_{\theta}^2\mathbb{H}_n(\theta^{\star}+u(\hat{\theta}_n-\theta^{\star}))D_n^{-1}du\right) D_n(\hat{\theta}_n-\theta^{\star}).
\end{align*}
We follow the classical route to show the asymptotic normality by verifying the following statements (for example, \cite[Lemma 4 and Lemma 5]{Kes97} or \cite[Lemma 2.2]{CHENG202473}):
\begin{align}
	&C_n(\theta^{\star}):= -D_n^{-1}\partial_{\theta}^2\mathbb{H}_n(\theta^{\star})D_n^{-1} \xrightarrow{p} -\mathcal{I}(\theta^{\star});
	\label{eq:an1}
	\\
	&\mathcal{I}_n(\theta^{\star}):=D_n^{-1}\partial_{\theta}\mathbb{H}_n(\theta^{\star}) \xrightarrow{\mathcal{L}} N(0,\mathcal{I}(\theta^{\star}));
	\label{eq:an2}
	\\
	&\sup_{\theta:|\theta|\leq \delta_n}|C_n(\theta^{\star}+\theta)-C_n(\theta^{\star})| \xrightarrow{p} 0 \text{ where $\delta_n \to 0$}.
	\label{eq:an3}
\end{align}
Recall the expression of $\mathbb{H}_n(\theta)$ in \eqref{eq:GQLF}, we list here the first and second order derivatives of $\mathbb{H}_n(\theta)$. For $i \in S$,
\begin{align*}
	&\partial_{\alpha_i}\mathbb{H}_n(\theta) = \sum_{j=1}^{n}\left(\frac{\left(X_{t_j}-\mu_{j-1}(\alpha)\right)\partial_{\alpha_i}b_{j-1}(\alpha)}{\sigma_{j-1}^2(\gamma)}\right)I_{\{Z_{t_{j-1}}=i\}},
	\\
	&\partial_{\gamma_i}\mathbb{H}_n(\theta) =  \frac{1}{2}\sum_{j=1}^{n} \left(\frac{\partial_{\gamma_i}\sigma_{j-1}^2(\gamma)\left(X_{t_j}-\mu_{j-1}(\alpha)\right)^2}{\sigma_{j-1}^4(\gamma)h} 
	-
	\frac{\partial_{\gamma_i}\sigma_{j-1}^2(\gamma)}{\sigma_{j-1}^2(\gamma)}
	\right)I_{\{Z_{t_{j-1}}=i\}},
	\\
	&\partial_{\alpha_i}^2\mathbb{H}_n(\theta) = \sum_{j=1}^{n}\left(
	-\frac{h\left(\partial_{\alpha_i}b_{j-1}(\alpha)\right)^2}{\sigma_{j-1}^2(\gamma)}
	+
	\frac{\left(X_{t_j}-\mu_{j-1}(\alpha)\right)\partial_{\alpha_i}^2b_{j-1}(\alpha)}{\sigma_{j-1}^2(\gamma)}
	\right)I_{\{Z_{t_{j-1}}=i\}},
	\\
	&\partial_{\gamma_i}^2\mathbb{H}_n(\theta) =  \frac{1}{2}\sum_{j=1}^{n} \Bigg(
	\frac{\left(\sigma_{j-1}^2(\gamma)\partial_{\gamma_i}^2\sigma_{j-1}^2(\gamma)-2\left(\partial_{\gamma_i}\sigma_{j-1}^2(\gamma)\right)^2\right)\left(X_{t_j}-\mu_{j-1}(\alpha)\right)^2}{\sigma_{j-1}^6(\gamma)h}
	\\
	&\quad\quad\quad\quad\quad\quad\quad\quad -
	\frac{\sigma_{j-1}^2(\gamma)\partial_{\gamma_i}^2\sigma_{j-1}^2(\gamma)-\left(\partial_{\gamma_i}\sigma_{j-1}^2(\gamma)\right)^2}{\sigma_{j-1}^4(\gamma)}
	\Bigg) I_{\{Z_{t_{j-1}}=i\}},
	\\
	&\partial_{\alpha_i}\partial_{\gamma_i}\mathbb{H}_n(\theta) = \sum_{j=1}^{n} \left(-\frac{2\partial_{\gamma_i}\sigma_{j-1}^2(\gamma)\left(X_{t_j}-\mu_{j-1}(\alpha)\right)\partial_{\alpha_i}b_{j-1}(\alpha)}{\sigma_{j-1}^4(\gamma)}
	\right)I_{\{Z_{t_{j-1}}=i\}}.
\end{align*}

We first prove 
\eqref{eq:an1}. 
For $i \in S$, it followed from Lemma \ref{lem:disergodic}, Lemma \ref{lem:estimates} and Assumption \ref{Ass:finite moment} that
\begin{align*}
	\frac{1}{nh}\partial_{\alpha_i}^2\mathbb{H}_n(\theta^{\star}) 
	&=  \frac{1}{n}\sum_{j=1}^{n}\left(
	-\frac{\left(\partial_{\alpha_i}b_{j-1}(\alpha^{\star})\right)^2}{\sigma_{j-1}^2(\gamma^{\star})}
	\right)I_{\{Z_{t_{j-1}}=i\}} + O_p(\sqrt{h})
	\\
	&= \sum_{k=1}^{N}\int_{\mathbb{R}}\left(-\frac{\left(\partial_{\alpha_i}b(x,k,\alpha_k^{\star})\right)^2}{\sigma^2(x,k,\gamma_k^{\star})}I_{\{k=i\}}\right) \nu_{\theta^{\star}}(dx,k) + o_p(1)
	\\
	&= -\int_{\mathbb{R}}\frac{\left(\partial_{\alpha_i}b(x,i,\alpha_i^{\star})\right)^2}{\sigma^2(x,i,\gamma_i^{\star})} \nu_{\theta^{\star}}(dx,i) + o_p(1),
\end{align*}
\begin{align*}
	\frac{1}{n}\partial_{\gamma_i}^2\mathbb{H}_n(\theta^{\star}) 
	&= \frac{1}{2n}\sum_{j=1}^{n} \Bigg(
	-\frac{\left(\partial_{\gamma_i}\sigma_{j-1}^2(\gamma^{\star})\right)^2}{\sigma_{j-1}^4(\gamma^{\star})}
	\Bigg) I_{\{Z_{t_{j-1}}=i\}}
	+ O_p(\sqrt{h})
	\\
	&= \sum_{k=1}^{N} \int_{\mathbb{R}}\left(-\frac{\left(\partial_{\gamma_i}\sigma^2(x,k,\gamma_k^{\star})\right)^2}{2\sigma^4(x,k,\gamma_i^{\star})}I_{\{k=i\}}\right) \nu_{\theta^{\star}}(dx,k) + o_p(1)
	\\
	&= -\int_{\mathbb{R}}\frac{\left(\partial_{\gamma_i}\sigma^2(x,i,\gamma_i^{\star})\right)^2}{2\sigma^4(x,i,\gamma_i^{\star})} \nu_{\theta^{\star}}(dx,i) + o_p(1),
\end{align*}
and
\begin{align*}
	\frac{1}{n\sqrt{h}}\partial_{\alpha_i}\partial_{\gamma_i}\mathbb{H}_n(\theta^{\star}) 
	&=O_p(h).
\end{align*}
The above convergences are enough to conclude \eqref{eq:an1}.

To prove 
\eqref{eq:an2}, we write
\begin{align*}
	\xi_{1,i,j}&=\frac{1}{\sqrt{nh}}\left(\frac{\left(X_{t_j}-\mu_{j-1}(\alpha^{\star})\right)\partial_{\alpha_i}b_{j-1}(\alpha^{\star})}{\sigma_{j-1}^2(\gamma^{\star})}\right)I_{\{Z_{t_{j-1}}=i\}},
	\\
	\xi_{2,i,j}&= \frac{1}{2\sqrt{n}}\left(\frac{\partial_{\gamma_i}\sigma_{j-1}^2(\gamma^{\star})\left(X_{t_j}-\mu_{j-1}(\alpha^{\star})\right)^2}{\sigma_{j-1}^4(\gamma^{\star})h} 
	-
	\frac{\partial_{\gamma_i}\sigma_{j-1}^2(\gamma^{\star})}{\sigma_{j-1}^2(\gamma^{\star})}
	\right)I_{\{Z_{t_{j-1}}=i\}},
\end{align*} 
where $i \in S$.
To establish the convergence indicated in \eqref{eq:an2}, we apply the central limit theorem to triangular arrays of random variables \cite[Lemma 3.6]{kessler2012statistical}.
Note that $\xi_{a,k,j} \xi_{b,l,j} = 0$ for all $a, b \in \{1,2\}$ and $k \neq l$. Thus,
\begin{equation}
	\sum_{j=1}^{n} \mathbb{E}_{\theta^{\star}}^{j-1}\left( \xi_{a,k,j}\xi_{b,l,j} \right)\xrightarrow{p} 0,
\end{equation} for $k \neq l$.
Therefore, \eqref{eq:an2} follows on showing that for $i \in S$,
\begin{align}
	&\sum_{j=1}^{n} \mathbb{E}_{\theta^{\star}}^{j-1}\left( \xi_{1,i,j} \right)\xrightarrow{p} 0, \quad
	\sum_{j=1}^{n} \mathbb{E}_{\theta^{\star}}^{j-1}\left( \xi_{2,i,j} \right)\xrightarrow{p} 0,
	\label{eq:an2c1}
	\\
	&\sum_{j=1}^{n} \mathbb{E}_{\theta^{\star}}^{j-1}\left( \xi_{1,i,j}^2 \right)\xrightarrow{p} G_{1,i}(\theta^{\star}), \quad
	\sum_{j=1}^{n} \mathbb{E}_{\theta^{\star}}^{j-1}\left( \xi_{2,i,j}^2 \right)\xrightarrow{p} G_{2,i}(\theta^{\star}),
	\notag
	\\
	&\sum_{j=1}^{n} \mathbb{E}_{\theta^{\star}}^{j-1}\left( \xi_{1,i,j}\xi_{2,i,j} \right)\xrightarrow{p} 0,
	\label{eq:an2c2}
	\\
	&\sum_{j=1}^{n} \mathbb{E}_{\theta^{\star}}^{j-1}\left( |\xi_{1,i,j}|^4 \right)\xrightarrow{p} 0, \quad
	\sum_{j=1}^{n} \mathbb{E}_{\theta^{\star}}^{j-1}\left( |\xi_{2,i,j}|^4 \right)\xrightarrow{p} 0.
	\label{eq:an2c3}
\end{align}
We first look at \eqref{eq:an2c1}. By Lemma \ref{lem:disergodic}, Lemma \ref{lem:estimates}, Assumption \ref{Ass:finite moment} and $nh^2 \to 0$, it is easy to see that 
\begin{align*}
	\sum_{j=1}^{n} \mathbb{E}_{\theta^{\star}}^{j-1}\left( \xi_{1,i,j} \right) &= \frac{1}{n} \sum_{j=1}^{n} \sqrt{n}h R_{j-1} = o_p(1),
	\\
	\sum_{j=1}^{n} \mathbb{E}_{\theta^{\star}}^{j-1}\left( \xi_{2,i,j} \right) &= \frac{1}{n} \sum_{j=1}^{n} \sqrt{n}h^{3/2}R_{j-1} = o_p(1).
\end{align*}
For the second part \eqref{eq:an2c2}. Similar computation shows that 
\begin{align*}
	\sum_{j=1}^{n} \mathbb{E}_{\theta^{\star}}^{j-1}\left( \xi_{1,i,j}^2 \right) &= \frac{1}{n} \sum_{j=1}^{n} \mathbb{E}_{\theta^{\star}}^{j-1}\left( \frac{1}{h}\left(\frac{\left(X_{t_j}-\mu_{j-1}(\alpha^{\star})\right)^2(\partial_{\alpha_i}b_{j-1}(\alpha^{\star}))^2}{\sigma_{j-1}^4(\gamma^{\star})}\right)I_{\{Z_{t_{j-1}}=i\}} \right)
	\\
	&= \frac{1}{n} \sum_{j=1}^{n} \left(\left(\frac{(\partial_{\alpha_i}b_{j-1}(\alpha^{\star}))^2}{\sigma_{j-1}^2(\gamma^{\star})}\right)I_{\{Z_{t_{j-1}}=i\}} \right) + O_p(h)
	\\
	&= G_{1,i}(\theta^{\star}) + o_p(1),
	\\
	\sum_{j=1}^{n} \mathbb{E}_{\theta^{\star}}^{j-1}\left( \xi_{2,i,j}^2 \right) &= \frac{1}{4n} \sum_{j=1}^{n} \mathbb{E}_{\theta^{\star}}^{j-1} \left(\frac{\partial_{\gamma_i}\sigma_{j-1}^2(\gamma^{\star})\left(X_{t_j}-\mu_{j-1}(\alpha^{\star})\right)^2}{\sigma_{j-1}^4(\gamma^{\star})h} 
	-
	\frac{\partial_{\gamma_i}\sigma_{j-1}^2(\gamma^{\star})}{\sigma_{j-1}^2(\gamma^{\star})}
	\right)^2I_{\{Z_{t_{j-1}}=i\}}
	\\
	&= \frac{1}{4n} \sum_{j=1}^{n} \mathbb{E}_{\theta^{\star}}^{j-1} \Bigg\{\frac{(\partial_{\gamma_i}\sigma_{j-1}^2(\gamma^{\star}))^2\left(X_{t_j}-\mu_{j-1}(\alpha^{\star})\right)^4}{\sigma_{j-1}^8(\gamma^{\star})h^2} 
	+
	\frac{(\partial_{\gamma_i}\sigma_{j-1}^2(\gamma^{\star}))^2}{\sigma_{j-1}^4(\gamma^{\star})}
	\\
	&\quad-
	2\frac{(\partial_{\gamma_i}\sigma_{j-1}^2(\gamma^{\star}))^2\left(X_{t_j}-\mu_{j-1}(\alpha^{\star})\right)^2}{\sigma_{j-1}^6(\gamma^{\star})h}
	\Bigg\}I_{\{Z_{t_{j-1}}=i\}}
	\\
	&= \frac{1}{n} \sum_{j=1}^{n} \mathbb{E}_{\theta^{\star}}^{j-1} \left(\left(
	\frac{(\partial_{\gamma_i}\sigma_{j-1}^2(\gamma^{\star}))^2}{2\sigma_{j-1}^4(\gamma^{\star})}
	\right)I_{\{Z_{t_{j-1}}=i\}}\right) + O_p(\sqrt{h})
	\\
	&= G_{2,i}(\theta^{\star}) + o_p(1),
	\end{align*} and
	\begin{align*}
	\sum_{j=1}^{n} \mathbb{E}_{\theta^{\star}}^{j-1}\left( \xi_{1,i,j}\xi_{2,i,j} \right) &= \frac{1}{2n\sqrt{h}} \sum_{j=1}^{n} \mathbb{E}_{\theta^{\star}}^{j-1} \Bigg\{
	\frac{\left(X_{t_j}-\mu_{j-1}(\alpha^{\star})\right)^3\partial_{\alpha_i}b_{j-1}(\alpha^{\star})\partial_{\gamma_i}\sigma_{j-1}^2(\gamma^{\star})}{\sigma_{j-1}^6(\gamma^{\star})h}
	\\
	&\quad -
	\frac{\left(X_{t_j}-\mu_{j-1}(\alpha^{\star})\right)\partial_{\alpha_i}b_{j-1}(\alpha^{\star})\partial_{\gamma_i}\sigma_{j-1}^2(\gamma^{\star})}{\sigma_{j-1}^4(\gamma^{\star})}
	\Bigg\}I_{\{Z_{t_{j-1}}=i\}}
	\\
	&= \frac{1}{2n} \sum_{j=1}^{n} h R_{j-1} = o_p(1).
\end{align*}
By similar arguments above, it is easy to show the third part \eqref{eq:an2c3} as following
\begin{align*}
	&\sum_{j=1}^{n} \mathbb{E}_{\theta^{\star}}^{j-1}\left( |\xi_{1,i,j}|^4 \right) = \frac{1}{n^2}\sum_{j=1}^{n} (R_{j-1}(\theta^{\star}) + hR_{j-1} = o_p(1),
	\\
	&\sum_{j=1}^{n} \mathbb{E}_{\theta^{\star}}^{j-1}\left( |\xi_{2,i,j}|^4 \right) = \frac{1}{n^2} \sum_{j=1}^{n} R_{j-1}= o_p(1).
\end{align*}

To show \eqref{eq:an3}, we observe that 
\begin{align}
	\sup_{\theta:|\theta|\leq \delta_n}|C_n(\theta^{\star}+\theta)-C_n(\theta^{\star})| \leq \delta_n \sup_{\theta} |\partial_{\theta}C_n(\theta)|.
\end{align} By applying Lemma \ref{lem:disergodic} and Lemma \ref{lem:estimates} to the partial derivative $\partial_{\theta}C_n(\theta)$, it not difficult to obtain that $\sup_{\theta} |\partial_{\theta}C_n(\theta)| = O_p(1)$.
The proof is complete.

\subsection{Proof of Corollary \ref{cor:consistencyI}}

The result follows from an application of Theorem \ref{thm:AN} and Lemma \ref{lem:disergodic}.
First, observe that for $i \in S$,
\begin{align*}
	\hat{G}_{1,i}^{(n)} &= \frac{1}{n} \sum_{j=1}^{n} \frac{\left(\partial_{\alpha_i}b(X_{t_{j-1}},Z_{t_{j-1}},\hat{\alpha}_{i,n})\right)^2}{\sigma^2(X_{t_{j-1}},Z_{t_{j-1}},\hat{\gamma}_{i,n})} I_{\{Z_{t_{j-1}}=i\}}
	\\
	&= \frac{1}{n} \sum_{j=1}^{n} \frac{\left(\partial_{\alpha_i}b(X_{t_{j-1}},Z_{t_{j-1}},\alpha_{Z_{t_{j-1}}})\right)^2}{\sigma^2(X_{t_{j-1}},Z_{t_{j-1}},\gamma_{Z_{t_{j-1}}})} I_{\{Z_{t_{j-1}}=i\}}
	\\
	&\quad + \Bigg(\frac{1}{n} \sum_{j=1}^{n} \frac{\left(\partial_{\alpha_i}b(X_{t_{j-1}},Z_{t_{j-1}},\hat{\alpha}_{i,n})\right)^2}{\sigma^2(X_{t_{j-1}},Z_{t_{j-1}},\hat{\gamma}_{i,n})} I_{\{Z_{t_{j-1}}=i\}}
	\\
	&\quad-
	\frac{1}{n} \sum_{j=1}^{n} \frac{\left(\partial_{\alpha_i}b(X_{t_{j-1}},Z_{t_{j-1}},\alpha_{Z_{t_{j-1}}})\right)^2}{\sigma^2(X_{t_{j-1}},Z_{t_{j-1}},\hat{\gamma}_{i,n})} I_{\{Z_{t_{j-1}}=i\}}\Bigg)
	\\
	&\quad +
	\Bigg(\frac{1}{n} \sum_{j=1}^{n} \frac{\left(\partial_{\alpha_i}b(X_{t_{j-1}},Z_{t_{j-1}},\alpha_{Z_{t_{j-1}}})\right)^2}{\sigma^2(X_{t_{j-1}},Z_{t_{j-1}},\hat{\gamma}_{i,n})} I_{\{Z_{t_{j-1}}=i\}}
	\nonumber\\
	&\quad -
	\frac{1}{n} \sum_{j=1}^{n} \frac{\left(\partial_{\alpha_i}b(X_{t_{j-1}},Z_{t_{j-1}},\alpha_{Z_{t_{j-1}}})\right)^2}{\sigma^2(X_{t_{j-1}},Z_{t_{j-1}},\gamma_{Z_{t_{j-1}}})} I_{\{Z_{t_{j-1}}=i\}}\Bigg)
	\\
	&=:
	\frac{1}{n} \sum_{j=1}^{n} \frac{\left(\partial_{\alpha_i}b(X_{t_{j-1}},Z_{t_{j-1}},\alpha_{Z_{t_{j-1}}})\right)^2}{\sigma^2(X_{t_{j-1}},Z_{t_{j-1}},\gamma_{Z_{t_{j-1}}})} I_{\{Z_{t_{j-1}}=i\}}
	+ S_{1,n} + S_{2,n}.
\end{align*}
Lemma \ref{lem:disergodic} shows that
\begin{equation*}
	\frac{1}{n} \sum_{j=1}^{n} \frac{\left(\partial_{\alpha_i}b(X_{t_{j-1}},Z_{t_{j-1}},\alpha_{Z_{t_{j-1}}})\right)^2}{\sigma^2(X_{t_{j-1}},Z_{t_{j-1}},\gamma_{Z_{t_{j-1}}})} I_{\{Z_{t_{j-1}}=i\}}
	\xrightarrow{p}
	G_{1,i}(\theta^{\star}).
\end{equation*}
Write $\psi(x,i,\alpha_i) = \left(\partial_{\alpha_i}b(x,i,\alpha_i)\right)^2$. By Assumption \ref{ass:smoothness}, the Taylor expansion, the moment assumption in Assumption \ref{Ass:finite moment}, and the consistency, we obtain
\begin{align*}
	&|S_{1,n}| 
	\\
	&\leq \frac{1}{n} \sum_{j=1}^{n} \left|\frac{\psi(X_{t_{j-1}},Z_{t_{j-1}},\hat{\alpha}_{i,n})-\psi(X_{t_{j-1}},Z_{t_{j-1}},\alpha_{Z_{t_{j-1}}})}{\sigma^2(X_{t_{j-1}},Z_{t_{j-1}},\hat{\gamma}_{i,n})} \right|I_{\{Z_{t_{j-1}}=i\}}
	\\
	&\leq \frac{1}{n} \sum_{j=1}^{n} \sigma^{-2}(X_{t_{j-1}},i,\hat{\gamma}_{i,n}) \int_{0}^{1}\left|\partial_{\alpha_i} \psi(X_{t_{j-1}},i,\alpha_i+u(\hat{\alpha}_{i,n}-\alpha_i)) \right|duI_{\{Z_{t_{j-1}}=i\}} \left|\hat{\alpha}_{i,n}-\alpha_i\right|
	\\
	&\leq  \frac{1}{n} \sum_{j=1}^{n} C(1+|X_{t_{j-1}}|) ^C\left|\hat{\alpha}_{i,n}-\alpha_i\right|
	\\
	&= O_p(1)o_p(1) = o_p(1).
\end{align*}
Using the same approach, we can also show $S_{2,n} = o_p(1)$. Consequently, we conclude that
\begin{equation}
	\hat{G}_{1,i}^{(n)} \xrightarrow{p}  G_{1,i}(\theta^{\star}).
	\label{eq:cor1c1}
\end{equation}
Following the above approach, it is straightforward to establish that
\begin{equation}
	\hat{G}_{2,i}^{(n)} \xrightarrow{p}  G_{2,i}(\theta^{\star})
	\label{eq:cor1c2}
\end{equation} for $i \in S$. The consistency $\hat{\mathcal{I}}_n \xrightarrow{p} \mathcal{I}(\theta^{\star})$ is thus followed from \eqref{eq:cor1c1} and \eqref{eq:cor1c2}.

\subsection{Proof of Proposition \ref{prop:ex1}}

First, we prove the following auxiliary lemma.

\begin{lemma} \label{lem:ZIDL}
	For any $i \in S$, we have
	\begin{align}
		\frac{1}{n}\sum_{j=1}^{n}I_{\{Z_{t_{j-1}}=i\}}  = \frac{1}{T}\int_{0}^{T}I_{\{Z_s=i\}}ds + o_p(1).
		\label{lem:ZIDL-1}
	\end{align}
\end{lemma}
\begin{proof}
	Observe that
	\begin{align*}
		\frac{1}{n}\sum_{j=1}^{n}I_{\{Z_{t_{j-1}}=i\}}-\frac{1}{nh}\int_{0}^{nh}I_{\{Z_s=i\}}ds
		=
		\frac{1}{nh}\sum_{j=1}^{n}\int_{t_{j-1}}^{t_j}\left(I_{\{Z_{t_{j-1}}=i\}}-I_{\{Z_s=i\}}\right)ds.
	\end{align*}
	By Equation (2.29) and \cite[Lemma A.5]{yin2009hybrid},  it follows that 
	\begin{equation}
		\mathbb{E}_{Q^{\star}}^{j-1} \left( I_{\{Z_{t_j}=i\}} \right) = I_{\{Z_{t_{j-1}}=i\}} + 
		I_{j-1}(i, Q^{\star}),
		\label{eq:ZcondiE1}
	\end{equation} where
	\begin{equation*}
		I_{j-1}(i, Q^{\star}) := \mathbb{E}_{Q^{\star}}^{j-1} \left(\int_{t_{j-1}}^{t_j} \left( \sum_{l=1}^{N} I_{\{Z_{s}=l\}} q^{\star}_{li} \right) ds\right)
	\end{equation*}
	and 
	\begin{equation}
		\max_i \sup_j I_{j-1}(i, Q^{\star}) = O_p(h).
		\label{eq:ZcondiE2}
	\end{equation}
	By applying \eqref{eq:ZcondiE1} and \eqref{eq:ZcondiE2}, 
	we have
	\begin{align*}
		& \frac{1}{nh}\sum_{j=1}^{n}\mathbb{E}_{Q^{\star}}^{j-1}\left(\int_{t_{j-1}}^{t_j}\left(I_{\{Z_{t_{j-1}}=i\}}-I_{\{Z_s=i\}}\right)ds\right)
		\nonumber\\
		&=
		\frac{1}{nh}\sum_{j=1}^{n}\int_{t_{j-1}}^{t_j}\mathbb{E}_{Q^{\star}}^{j-1}\left(I_{\{Z_{t_{j-1}}=i\}}-I_{\{Z_s=i\}}\right)ds
		\\
		&=\frac{1}{nh}\sum_{j=1}^{n}\int_{t_{j-1}}^{t_j}\left(I_{\{Z_{t_{j-1}}=i\}}-(I_{\{Z_{t_{j-1}}=i\}} + I_{j-1}(i, Q^{\star}))\right)ds
		\\
		&=\frac{1}{nh}\sum_{j=1}^{n}\int_{t_{j-1}}^{t_j} I_{j-1}(i, Q^{\star})ds
		= O_p(h).
	\end{align*}
	By Jensen's inequality, we also have
	\begin{align*}
		&\frac{1}{n^2h^2}\sum_{j=1}^{n}\mathbb{E}_{Q^{\star}}^{j-1}\left(\int_{t_{j-1}}^{t_j}\left(I_{\{Z_{t_{j-1}}=i\}}-I_{\{Z_s=i\}}\right)ds\right)^2
		\\
		&\leq \frac{1}{n^2h}\sum_{j=1}^{n}\mathbb{E}_{Q^{\star}}^{j-1}\int_{t_{j-1}}^{t_j}\left(I_{\{Z_{t_{j-1}}=i\}}-I_{\{Z_s=i\}}\right)^2ds
		\\
		&= O_p(1/n).
	\end{align*}
	Thus we obtain \eqref{lem:ZIDL-1}.
\end{proof}

Now we prove Proposition \ref{prop:ex1}. Note that
\begin{align}
	\sqrt{nh}(\hat{\pi}_k-\pi_k) = \sqrt{nh}\left(\hat{\pi}_k-\frac{1}{T}\int_{0}^{T}I_{\{Z_s=k\}}ds\right) + \sqrt{nh}\left(\frac{1}{T}\int_{0}^{T}I_{\{Z_s=k\}}ds - \pi_k\right).
	\label{eq:prop:1}
\end{align}
By similar computations in Lemma \ref{lem:ZIDL}, it is easy to obtain
\begin{align*}
	&\sum_{j=1}^{n}\mathbb{E}_{Q^{\star}}^{j-1}\left(\frac{1}{\sqrt{nh}}\int_{t_{j-1}}^{t_j}\left(I_{\{Z_{t_{j-1}}=k\}}-I_{\{Z_s=k\}}\right)ds\right) = O_p(\sqrt{nh^3})=o_p(1),
	\\
	&\sum_{j=1}^{n}\mathbb{E}_{Q^{\star}}^{j-1}\left(\frac{1}{\sqrt{nh}}\int_{t_{j-1}}^{t_j}\left(I_{\{Z_{t_{j-1}}=k\}}-I_{\{Z_s=k\}}\right)ds\right)^2 = O_p(h)=o_p(1),
\end{align*}
hence
\begin{equation}
	\sqrt{nh}\left(\hat{\pi}_k-\frac{1}{T}\int_{0}^{T}I_{\{Z_s=k\}}ds\right) \xrightarrow{p} 0.
	\label{eq:prop:2}
\end{equation}
For the second term on the right-hand side of \eqref{eq:prop:1}, we apply a central limit theorem for exponentially ergodic continuous-time Markov chains (see, for example, \cite[Theorem 3.1]{liu2015central}). 
It can be used by verifying
\begin{equation*}
	\sum_{i\in S} |I_{\{i=k\}}|^4 \pi_i = \sum_{i\in S} I_{\{i=k\}} \pi_i = \pi_k < \infty.
\end{equation*}
Recall that $\tau_{i}$ is the first return time to state $i$. So for any fixed $i_0 \in S$, \cite[Theorem 3.1]{liu2015central} ensures
\begin{align} \label{eq:prop:3}
	\sqrt{nh}\left(\frac{1}{T}\int_{0}^{T}I_{\{Z_s=k\}}ds - \pi_k\right)\xrightarrow{\mathcal{L}} N(0,V_{\pi_k}),
\end{align}
where
\begin{align*}
	&V_{\pi_k}= 2 \sum_{i \in \{1,2\}} (I_{\{i=k\}}-\pi_k)F_i\pi_i ,
	\\
	&F_i = \mathbb{E}_{Q^{\star}}\left(\int_{0}^{\tau_{i_0}}(I_{\{Z_s=k\}}-\pi_k)ds \,\middle|\, Z_0 = i\right).
\end{align*}
Substituting \eqref{eq:prop:2} and \eqref{eq:prop:3} in \eqref{eq:prop:1}, we obtain the desired result.

\subsection{Proof of Theorem \ref{thm:consisQ}}

Given the relations $\hat{q}^{(n)}_{ii} = -\sum_{i \neq k} \hat{q}^{(n)}_{ik}$ and $q^{\star}_{ii} = -\sum_{i \neq k} q^{\star}_{ik}$, it is sufficient to show $\hat{q}^{(n)}_{ik} \xrightarrow{p} q^{\star}_{ik}$ for $i \neq k$. From equation \eqref{eq:estQn}, we have for $i \neq k$
\begin{align*}
	\hat{q}^{(n)}_{ik} = \frac{K_{ik}(n)}{h \sum_{l=1}^{N} K_{il}(n)} = \frac{\frac{1}{nh} K_{ik}(n)}{\frac{1}{n} \sum_{l=1}^{N} K_{il}(n)}.
\end{align*}

We first consider
\begin{align*}
	\frac{1}{n} K_{ii}(n) = \frac{1}{n} \sum_{j=1}^{n} I_{\{Z_{t_{j-1}}=i\}} I_{\{Z_{t_j}=i\}},\qquad i\in S.
\end{align*}
Additionally, recall that under Assumption \ref{Ass:finite moment}, $Z_t$ is ergodic, i.e., there is a unique stationary distribution $\pi = (\pi_1, \ldots, \pi_N)$ of $Z_t$. Lemma \ref{lem:ZIDL} and the ergodic theorem for $Z_t$ (see, for example, \cite[Chapter 9.5]{grimmett2020probability} or  \cite[Theorem 3.8.1]{norris1998markov}) state that 
\begin{equation}
	\frac{1}{n} \sum_{j=1}^{n} I_{\{Z_{t_{j-1}}=i\}} \xrightarrow{p} \sum_{k=1}^{N} I_{\{k=i\}} \pi_k.
	\label{eq:MCergodic}
\end{equation} Thus, by \eqref{eq:ZcondiE1}, \eqref{eq:ZcondiE2} and the above ergodic theorem \eqref{eq:MCergodic}, we obtain
\begin{align*}
	\frac{1}{n} \sum_{j=1}^{n} \mathbb{E}_{Q^{\star}}^{j-1} \left( I_{\{Z_{t_{j-1}}=i\}} I_{\{Z_{t_j}=i\}} \right)
	&= \frac{1}{n} \sum_{j=1}^{n} I_{\{Z_{t_{j-1}}=i\}} \mathbb{E}_{Q^{\star}}^{j-1} \left( I_{\{Z_{t_j}=i\}} \right) \\
	&= \frac{1}{n} \sum_{j=1}^{n} I_{\{Z_{t_{j-1}}=i\}} (I_{\{Z_{t_{j-1}}=i\}} + I_{j-1}(i, Q^{\star})) \\
	&= \frac{1}{n} \sum_{j=1}^{n} I_{\{Z_{t_{j-1}}=i\}} + O_p(h) \\
	&= \sum_{k=1}^{N} I_{\{i=k\}} \pi_k + o_p(1) \\
	&= \pi_i + o_p(1),
\end{align*}
and
\begin{align*}
	\sum_{j=1}^{n} \frac{1}{n^2} \mathbb{E}_{Q^{\star}}^{j-1} \left( I_{\{Z_{t_{j-1}}=i\}} I_{\{Z_{t_j}=i\}} \right)^2
	&= \frac{1}{n^2} \sum_{j=1}^{n} \mathbb{E}_{Q^{\star}}^{j-1} \left( I_{\{Z_{t_{j-1}}=i\}} I_{\{Z_{t_j}=i\}} \right) \\
	&= O_p \left( \frac{1}{n} \right).
\end{align*}
Then \cite[Lemma 9]{GenonCatalot1993OnTE} implies 
\begin{equation}
	\frac{1}{n} K_{ii}(n) \xrightarrow{p} \pi_i.
	\label{eq:thmQc1}
\end{equation}

Now we look at
\begin{equation*}
	\frac{1}{nh} K_{ik}(n) = \frac{1}{nh} \sum_{j=1}^{n} I_{\{Z_{t_{j-1}}=i\}} I_{\{Z_{t_j}=k\}},\qquad i \neq k.
\end{equation*}
Note that by \eqref{eq:ZcondiE1} and \eqref{eq:ZcondiE2},
\begin{align*}
	\frac{1}{n^2 h^2} \sum_{j=1}^{n} \mathbb{E}_{Q^{\star}}^{j-1} \left( I_{\{Z_{t_{j-1}}=i\}} I_{\{Z_{t_j}=k\}} \right)^2
	&= \frac{1}{nh^2} \frac{1}{n} \sum_{j=1}^{n} \mathbb{E}_{Q^{\star}}^{j-1} \left( I_{\{Z_{t_{j-1}}=i\}} I_{\{Z_{t_j}=k\}} \right) \\
	&= O_p \left( \frac{1}{nh} \right).
\end{align*}
By \cite[Lemma A.5]{yin2009hybrid}, we have
\begin{equation*}
	\mathbb{E}_{Q^{\star}}^{j-1} \left( I_{\{Z_{t_j}=i\}} - I_{\{Z_{t_{j-1}}=i\}} - \int_{t_{j-1}}^{t_j} \left( \sum_{l=1}^{N} I_{\{Z_{s}=l\}} q^{\star}_{li} \right) ds \right) = 0
\end{equation*}
for all $i \in S$. Then, by the ergodic theorem, we have for $i\ne k$,
\begin{align*}
	& \frac{1}{nh} \sum_{j=1}^{n} \mathbb{E}_{Q^{\star}}^{j-1} \left( I_{\{Z_{t_{j-1}}=i\}} I_{\{Z_{t_j}=k\}} \right)
	\nonumber\\
	&= \frac{1}{nh} \sum_{j=1}^{n} I_{\{Z_{t_{j-1}}=i\}} \mathbb{E}_{Q^{\star}}^{j-1} \left( I_{\{Z_{t_j}=k\}} \right) \\
	&= \frac{1}{nh} \sum_{j=1}^{n} I_{\{Z_{t_{j-1}}=i\}} \mathbb{E}_{Q^{\star}}^{j-1} \left( I_{\{Z_{t_{j-1}}=k\}} + \int_{t_{j-1}}^{t_j} \left( \sum_{l=1}^{N} I_{\{Z_{s}=l\}} q^{\star}_{lk} \right) ds \right) \\
	&= \frac{1}{n} \sum_{j=1}^{n} I_{\{Z_{t_{j-1}}=i\}} \mathbb{E}_{Q^{\star}}^{j-1} \left( \frac{1}{h} \int_{t_{j-1}}^{t_j} \left( \sum_{l=1}^{N} I_{\{Z_{s}=l\}} q^{\star}_{lk} \right) ds \right)\\
	&= \frac{1}{n} \sum_{j=1}^{n} I_{\{Z_{t_{j-1}}=i\}} \left( \sum_{l=1}^{N} I_{\{Z_{t_{j-1}}=l\}} q^{\star}_{lk}\right)
	\\
	&\quad + \frac{1}{n} \sum_{j=1}^{n} I_{\{Z_{t_{j-1}}=i\}} \mathbb{E}_{Q^{\star}}^{j-1} \left( \frac{1}{h} \int_{t_{j-1}}^{t_j} \left( \sum_{l=1}^{N} I_{\{Z_{s}=l\}} q^{\star}_{lk} - \sum_{l=1}^{N} I_{\{Z_{t_{j-1}}=l\}} q^{\star}_{lk} \right) ds \right)\\
	&= \frac{1}{n} \sum_{j=1}^{n} I_{\{Z_{t_{j-1}}=i\}}  q^{\star}_{ik} + o_p(1) \\
	&= q^{\star}_{ik} \pi_i + o_p(1).
\end{align*}
Here, we used \eqref{eq:ZcondiE2} to get the second last identity.
Then \cite[Lemma 9]{GenonCatalot1993OnTE} implies
\begin{equation}
	\frac{1}{nh} K_{ik}(n) \xrightarrow{p} q^{\star}_{ik} \pi_i
	\label{eq:thmQc2}
\end{equation}
for $i \neq k$. From \eqref{eq:thmQc2}, we also have
\begin{equation}
	\frac{1}{n} K_{ik}(n) \xrightarrow{p} 0
	\label{eq:thmQc3}
\end{equation}
for $i \neq k$. Combining \eqref{eq:thmQc1}, \eqref{eq:thmQc2}, and \eqref{eq:thmQc3}, we obtain
\begin{equation*}
	\hat{q}^{(n)}_{ik} = \frac{\frac{1}{nh} K_{ik}(n)}{\frac{1}{n} \sum_{l=1}^{N} K_{il}(n)} \xrightarrow{p} q^{\star}_{ik}.
\end{equation*}

\section*{Acknowledgments}
This work was partially supported by WISE program (MEXT) at Kyushu University (YC), and by JST CREST Grant Number JPMJCR2115, Japan, and by JSPS KAKENHI Grant Number 23K22410 (HM).









\medskip
Received xxxx 20xx; revised xxxx 20xx; early access xxxx 20xx.
\medskip

\clearpage

\appendix   

\section*{Correction to Lemma 6.2}

This note corrects an error in Lemma 6.2 (p.18) of the paper. 
We keep the notation of the original article.

Concerned with \textbf{Lemma 6.2 (p.18)}, we need the two modifications given below.
With the corrected remainder order ``$h^2$'' in \eqref{hm:correction-1} below, the second convergence of Eq.(30) in p.24 should read
\begin{equation*}
	\sum_{j=1}^{n} \mathbb{E}_{\theta^{\star}}^{j-1}\left( \xi_{2,i,j} \right) = \frac{1}{n} \sum_{j=1}^{n} \sqrt{n}\,h R_{j-1} = o_p(1),
\end{equation*}
and all subsequent results, proofs, and conclusions remain unchanged.

\bigskip

\begin{enumerate}
	\item The first identity in Lemma 6.2 should be split into the two cases:
	\begin{align}
		&\mathbb{E}_{\theta^{\star}}^{j-1}\left[(X_{t_j}-\mu_{j-1}(\alpha))^2\right] = \sigma^2(X_{t_{j-1}},Z_{t_{j-1}},\gamma_{Z_{t_{j-1}}}^{\star})h + h^{3/2} R_{j-1}(\theta),
		\nonumber
		\\
		&\mathbb{E}_{\theta^{\star}}^{j-1}\left[(X_{t_j}-\mu_{j-1}(\alpha^{\star}))^2\right] = \sigma^2(X_{t_{j-1}},Z_{t_{j-1}},\gamma_{Z_{t_{j-1}}}^{\star})h + h^{2} R_{j-1}.
		\label{hm:correction-1}
	\end{align}
	
	\item In the proof (lines –3 on p.19), it was written that
	\begin{align*}
		&\mathbb{E}_{\theta^{\star}}^{j-1} \left( \int_{t_{j-1}}^{t_j} B_s ds \int_{t_{j-1}}^{t_j} \sigma(X_{t_{j-1}},Z_{t_{j-1}},\gamma_{Z_{t_{j-1}}}^{\star}) dw_s \right) = h^{3/2} R_{j-1}.
	\end{align*}
	This conditional expectation should read
	\begin{align}
		&\mathbb{E}_{\theta^{\star}}^{j-1} \left( \int_{t_{j-1}}^{t_j} B_s dw_s \int_{t_{j-1}}^{t_j} \sigma(X_{t_{j-1}},Z_{t_{j-1}},\gamma_{Z_{t_{j-1}}}^{\star}) dw_s \right) = h^{2} R_{j-1}.
		\label{hm:correction-3}
	\end{align}
	
\end{enumerate}
\begin{proof}[Proof of \eqref{hm:correction-3}]
	Write $B_s = B_{1,s} + B_{2,s}$, where
	\begin{align*}
		&B_{1,s} = \sigma(X_{s},Z_{s},\gamma_{Z_{s}}^{\star})-\sigma(X_{s},Z_{t_{j-1}},\gamma_{Z_{t_{j-1}}}^{\star}),
		\\
		&B_{2,s} = \sigma(X_{s},Z_{t_{j-1}},\gamma_{Z_{t_{j-1}}}^{\star})-\sigma(X_{t_{j-1}},Z_{t_{j-1}},\gamma_{Z_{t_{j-1}}}^{\star}).
	\end{align*}
	By a similar argument as in Lemma 6.2, we have 
	\begin{align*}
		|B_{1,s}| &= \left|\left(\sigma(X_{s},Z_{s},\gamma_{Z_{s}}^{\star})-\sigma(X_{s},Z_{t_{j-1}},\gamma_{Z_{t_{j-1}}}^{\star})\right)I_{\{Z_s \neq Z_{t_{j-1}}\}}\right|
		\\&\leq C(1+|X_{t_{j-1}}|+ |X_s - X_{t_{j-1}}|)^C I_{\{Z_s \neq Z_{t_{j-1}}\}}
		\\&\leq |R_{j-1}| I_{\{Z_s \neq Z_{t_{j-1}}\}} + C|X_s - X_{t_{j-1}}|^C I_{\{Z_s \neq Z_{t_{j-1}}\}}.
	\end{align*} Noting that $|\mathbb{E}_{\theta^{\star}}^{j-1}\left(R_{j-1} I_{\{Z_s \neq Z_{t_{j-1}}\}}\right)| \leq h |R_{j-1}|$ and using standard moment estimates gives
	\begin{equation*}
		\mathbb{E}_{\theta^{\star}}^{j-1}[B_{1,s}] = h R_{j-1}.
	\end{equation*}
	Next, a second‐order Taylor expansion gives 
	\begin{align*}
		B_{2,s} = \partial_{x}\sigma(X_{t_{j-1}},Z_{t_{j-1}},\gamma_{Z_{t_{j-1}}}^{\star})(X_s-X_{t_{j-1}}) +  M_{s,t_{j-1}} (X_s-X_{t_{j-1}})^2,
	\end{align*} where $M_{s,t_{j-1}}$ satisfies $|M_{s,t_{j-1}}| \leq C(1 + |X_{t_{j-1}}| + |X_s|)^C$ for some constant $C>0$.
	Taking conditional expectation and using standard moment estimates yields $|\mathbb{E}_{\theta^{\star}}^{j-1}[B_{2,s}]| \leq h |R_{j-1}|.$
	Hence, 
	\begin{equation*}
		\mathbb{E}_{\theta^{\star}}^{j-1}[B_{s}] = h R_{j-1}.
	\end{equation*}
	Finally, using It\^{o}'s isometry, we observe that
	\begin{align}
		& \mathbb{E}_{\theta^{\star}}^{j-1} \left( \int_{t_{j-1}}^{t_j} B_s dw_s \int_{t_{j-1}}^{t_j} \sigma(X_{t_{j-1}},Z_{t_{j-1}},\gamma_{Z_{t_{j-1}}}^{\star}) dw_s \right) 
		\nonumber\\
		&= \sigma(X_{t_{j-1}},Z_{t_{j-1}},\gamma_{Z_{t_{j-1}}}^{\star}) 
		\mathbb{E}_{\theta^{\star}}^{j-1} \left( \int_{t_{j-1}}^{t_j} B_s ds\right) 
		\nonumber\\
		&= \sigma(X_{t_{j-1}},Z_{t_{j-1}},\gamma_{Z_{t_{j-1}}}^{\star}) h^2 R_{j-1}
		= h^2 R_{j-1},
		\nonumber
	\end{align}
	hence \eqref{hm:correction-3}.
\end{proof}

\begin{proof}[Proof of \eqref{hm:correction-1}]
	Let $A_s(\theta^{\star}) := b(X_s, Z_s, \alpha_{Z_s}^{\star}) - b(X_{t_{j-1}}, Z_{t_{j-1}}, \alpha_{Z_{t_{j-1}}}^{\star})$.
	It suffices to show 
	\begin{equation}
		\mathbb{E}_{\theta^{\star}}^{j-1} \left( \int_{t_{j-1}}^{t_j} A_s(\theta^{\star}) ds \int_{t_{j-1}}^{t_j} \sigma(X_{t_{j-1}},Z_{t_{j-1}},\gamma_{Z_{t_{j-1}}}^{\star}) dw_s \right) = h^{2} R_{j-1}.
		\label{eq:A}
	\end{equation}
	By the same arguments as in Lemma 6.2, one has
	\begin{equation*}
		\mathbb{E}_{\theta^{\star}}^{j-1}[A_s(\theta^{\star})^2] = h R_{j-1}.
	\end{equation*} 
	Then, applying H\"{o}lder's inequality immediately yields the desired result \eqref{eq:A}.
\end{proof}


\begin{thebibliography}{99}
	
	
\bibitem{wreo25274} 
\newblock H. Alkhezi,
\newblock \emph{Bayesian Inference for Animal Movement in Continuous Time},
\newblock Ph.D thesis, University of Sheffield, 2019.

\bibitem{bao2016approximation} 
\newblock J. Bao, J. Shao and C. Yuan,
\newblock \doititle{Approximation of invariant measures for regime-switching diffusions},
\newblock \emph{Potential Analysis}, \textbf{44} (2016), 707--727.

\bibitem{BLACKWELL199787} 
\newblock P. Blackwell,
\newblock \doititle{Random diffusion models for animal movement},
\newblock \emph{Ecological Modelling}, \textbf{100} (1997), 87--102.

\bibitem{Blackwell2003} 
\newblock P. G. Blackwell,
\newblock \doititle{Bayesian Inference for Markov Processes with Diffusion and Discrete Components},
\newblock \emph{Biometrika}, \textbf{90} (2003), 613--627.

\bibitem{bladt2005statistical} 
\newblock M. Bladt and M. S\o rensen,
\newblock \doititle{Statistical inference for discretely observed Markov jump processes},
\newblock \emph{Journal of the Royal Statistical Society Series B: Statistical Methodology}, \textbf{67} (2005), 395--410.

\bibitem{CHENG202473} 
\newblock Y. Cheng, N. Hufnagel and H. Masuda,
\newblock \doititle{Estimation of ergodic square-root diffusion under high-frequency sampling},
\newblock \emph{Econometrics and Statistics}, \textbf{32} (2024), 73--87.

\bibitem{dos2018robust} 
\newblock G. Dos Reis and G. Smith,
\newblock \doititle{Robust and consistent estimation of generators in credit risk},
\newblock \emph{Quantitative Finance}, \textbf{18} (2018), 983--1001.

\bibitem{GenonCatalot1993OnTE} 
\newblock V. Genon-Catalot and J. Jacod,
\newblock \doititle{On the estimation of the diffusion coefficient for multi-dimensional diffusion processes},
\newblock \emph{Annales De L'Institut Henri Poincar\'e (Probabilit\'es et Statistiques)}, \textbf{29} (1993), 119--151.

\bibitem{grimmett2020probability} 
\newblock G. Grimmett and D. Stirzaker,
\newblock \emph{Probability and Random Processes},
\newblock Oxford University Press, 2020.

\bibitem{Hamiltion1989} 
\newblock J. D. Hamilton,
\newblock \doititle{A new approach to the economic analysis of nonstationary time series and the business cycle},
\newblock \emph{Econometrica}, \textbf{57} (1989), 357--384.

\bibitem{harris2013flexible} 
\newblock K. J. Harris and P. G. Blackwell,
\newblock \doititle{Flexible continuous-time modelling for heterogeneous animal movement},
\newblock \emph{Ecological Modelling}, \textbf{255} (2013), 29--37.

\bibitem{spuRs_book} 
\newblock O. Jones, R. Maillardet and A. Robinson,
\newblock \emph{Introduction to Scientific Programming and Simulation Using R},
\newblock 2\textsuperscript{nd} edition, CRC Press, Boca Raton, 2014.

\bibitem{Kes97} 
\newblock M. Kessler,
\newblock \doititle{Estimation of an ergodic diffusion from discrete observations},
\newblock \emph{Scandinavian Journal of Statistics}, \textbf{24} (1997), 211--229.

\bibitem{kessler2012statistical} 
\newblock M. Kessler, A. Lindner and M. S\o rensen,
\newblock \doititle{Statistical methods for stochastic differential equations},
\newblock \emph{Monographs on Statistics and Applied Probability}, \textbf{124} (2012), 7--12.

\bibitem{kunita1986tightness} 
\newblock H. Kunita,
\newblock \doititle{Tightness of probability measures in $D([0,t];C)$ and $D([0,t];D)$},
\newblock \emph{Journal of the Mathematical Society of Japan}, \textbf{38} (1986), 309--334.

\bibitem{liu2015central} 
\newblock Y. Liu and Y. Zhang,
\newblock \doititle{Central limit theorems for ergodic continuous-time Markov chains with applications to single birth processes},
\newblock \emph{Frontiers of Mathematics in China}, \textbf{10} (2015), 933--947.

\bibitem{mao2006stochastic} 
\newblock X. Mao and C. Yuan,
\newblock \emph{Stochastic Differential Equations with Markovian Switching},
\newblock Imperial College Press, 2006.

\bibitem{metzner2007generator} 
\newblock P. Metzner, E. Dittmer, T. Jahnke and C. Sch\"utte,
\newblock \doititle{Generator estimation of Markov jump processes},
\newblock \emph{Journal of Computational Physics}, \textbf{227} (2007), 353--375.

\bibitem{michelot2019} 
\newblock T. Michelot and P. G. Blackwell,
\newblock \doititle{State-switching continuous-time correlated random walks},
\newblock \emph{Methods in Ecology and Evolution}, \textbf{10} (2019), 637--649.

\bibitem{michelot2021varying} 
\newblock T. Michelot, R. Glennie, C. Harris and L. Thomas,
\newblock \doititle{Varying-coefficient stochastic differential equations with applications in ecology},
\newblock \emph{Journal of Agricultural, Biological and Environmental Statistics}, \textbf{26} (2021), 446--463.

\bibitem{norris1998markov} 
\newblock J. R. Norris,
\newblock \emph{Markov Chains},
\newblock Cambridge University Press, 1998.

\bibitem{patterson2017statistical} 
\newblock T. A. Patterson, A. Parton, R. Langrock, P. G. Blackwell, L. Thomas and R. King,
\newblock \doititle{Statistical modelling of individual animal movement: an overview of key methods and a discussion of practical challenges},
\newblock \emph{AStA Advances in Statistical Analysis}, \textbf{101} (2017), 399--438.

\bibitem{shao2018Euler} 
\newblock J. Shao,
\newblock \doititle{Invariant measures and Euler--Maruyama's approximations of state-dependent regime-switching diffusions},
\newblock \emph{SIAM Journal on Control and Optimization}, \textbf{56} (2018), 3215--3238.

\bibitem{van2000} 
\newblock A. W. van der Vaart,
\newblock \emph{Asymptotic Statistics},
\newblock (Cambridge Series in Statistical and Probabilistic Mathematics, vol.~3), Cambridge University Press, Cambridge, 1998.

\bibitem{yin2005numerical} 
\newblock G. G. Yin, X. Mao and K. Yin,
\newblock \doititle{Numerical approximation of invariant measures for hybrid diffusion systems},
\newblock \emph{IEEE Transactions on Automatic Control}, \textbf{50} (2005), 934--946.

\bibitem{yin2009hybrid} 
\newblock G. G. Yin and C. Zhu,
\newblock \emph{Hybrid Switching Diffusions: Properties and Applications},
\newblock (Stochastic Modelling and Applied Probability, vol.~63), Springer, 2009.

\bibitem{YUAN2003277} 
\newblock C. Yuan and X. Mao,
\newblock \doititle{Asymptotic stability in distribution of stochastic differential equations with Markovian switching},
\newblock \emph{Stochastic Processes and their Applications}, \textbf{103} (2003), 277--291.

\bibitem{ZHANG2022273} 
\newblock X. Zhang and T. Zhang,
\newblock \doititle{Barrier option pricing under a Markov regime switching diffusion model},
\newblock \emph{The Quarterly Review of Economics and Finance}, \textbf{86} (2022), 273--280.

\bibitem{Zhang2017} 
\newblock Z. Zhang and W. Wang,
\newblock \doititle{The stationary distribution of Ornstein--Uhlenbeck process with a two-state Markov switching},
\newblock \emph{Communications in Statistics -- Simulation and Computation}, \textbf{46} (2017), 4783--4794.

\bibitem{Yuhang2023} 
\newblock Y. Zhen and F. Xi,
\newblock \doititle{Least squares estimators for stochastic differential equations with Markovian switching},
\newblock \emph{Discrete and Continuous Dynamical Systems -- B}, \textbf{28} (2023), 4068--4086.


\end{thebibliography}
\end{document}